\documentclass[oneside,11pt]{article}

\usepackage[utf8]{inputenc}
\usepackage[T1]{fontenc}
\usepackage{amsthm}
\usepackage{amsmath}
\usepackage{amsfonts}
\usepackage{amssymb}
\usepackage{bbm}
\usepackage{graphicx} 
\usepackage{float}
\usepackage{booktabs}
\usepackage{multirow}
\usepackage{rotating}
\usepackage{array}
\usepackage{xcolor}
\usepackage[ruled]{algorithm2e}
\usepackage[labelformat=simple]{subcaption}

\DeclareCaptionLabelFormat{subcaptionlabel}{\normalfont(\textbf{#2}\normalfont)}
\newcolumntype{R}[1]{>{\raggedleft\let\newline\\\arraybackslash\hspace{0pt}}m{#1}}
\newcolumntype{L}[1]{>{\raggedright\let\newline\\\arraybackslash\hspace{0pt}}m{#1}}
\usepackage{breakcites}
\usepackage{nowidow}
\usepackage{enumitem}
\usepackage[top=1.5cm,bottom=2cm,right=2.5cm,left=2.5cm]{geometry}
\usepackage{titling}
\usepackage{natbib}
\usepackage{ifplatform}
\usepackage{textcomp}
\ifwindows
\fi
\allowdisplaybreaks

\newcommand{\lp}{\left(}
\newcommand{\rp}{\right)}
\newcommand{\lc}{\left[}
\newcommand{\rc}{\right]}

\newcommand{\rmd}{\mathrm{d}}
\newcommand{\R}{\mathbb{R}}
\newcommand{\bbL}{\mathbb{L}}
\newcommand{\cF}{\mathcal{F}}
\newcommand{\cD}{\mathcal{D}}
\newcommand{\cC}{\mathcal{C}}
\newcommand{\cR}{\mathcal{R}}
\newcommand{\cB}{\mathcal{B}}
\newcommand{\cN}{\mathcal{N}}
\newcommand{\wt}[1]{\widetilde{#1}}
\newcommand{\lrp}[1]{\left(#1\right)}
\newcommand{\lrc}[1]{\left[#1\right]}
\newcommand{\lrb}[1]{\left\{#1\right\}}
\newcommand{\Prob}[1]{\mathbb{P}\lp #1\rp}
\newcommand{\Pro}[2]{\mathbb{P}_{#2}\lp #1\rp}
\renewcommand{\Pr}{\mathbb{P}}
\newcommand{\Esp}[1]{\mathbb{E}\lc #1\rc}
\newcommand{\Es}[2]{\mathbb{E}_{#2}\lc #1\rc}
\newcommand{\E}{\mathbb{E}}

\newcommand{\Vs}[2]{\mathbb{V}\mathrm{ar}_{#2}\lc #1\rc}
\newtheorem{theorem}{Theorem}

\newtheorem{remark}{Remark}
\newtheorem{proposition}{Proposition}
\newtheorem{lemma}{Lemma}

\usepackage[pdftex,final,bookmarksnumbered,bookmarksopen=false,breaklinks,colorlinks]{hyperref}
\hypersetup{
	final,
	bookmarksnumbered=true,
	bookmarksopen=true,
	bookmarksopenlevel=0,
	unicode=false,
	pdftoolbar=true,
	pdfmenubar=true,
	pdffitwindow=false,
	pdftitle={Optimal stopping of an Ornstein–Uhlenbeck bridge},
	pdfdisplaydoctitle=true,
	pdftoolbar=true,
	pdfmenubar=true,
	pdflang={English},
	pdfauthor={Abel Azze, Bernardo D'Auria, Eduardo Garcia-Portugues},
	pdfsubject={arXiv paper},
	pdfcreator={Abel Azze},
	pdfproducer={Abel Azze},
	pdfkeywords={Free-boundary problem}{Optimal stopping}{Ornstein--Uhlenbeck bridge}{Time-inhomogeneity},
	pdfnewwindow=true,
	breaklinks=true,
	hidelinks,       
	linkcolor=black,
	citecolor=black,
	filecolor=magenta,
	urlcolor=cyan
}

\newif\ifmain
\maintrue
\newif\ifsupplement
\supplementtrue

\newif\iffigstabs
\figstabstrue

\begin{document}

\ifmain

%-----------------------------------------------%
\title{Optimal stopping of an Ornstein--Uhlenbeck bridge}
\setlength{\droptitle}{-1cm}
\predate{}%
\postdate{}%
\date{}
%-----------------------------------------------%

%-----------------------------------------------%
\author{Abel Azze$^{1,4}$, Bernardo D'Auria$^{2}$, and Eduardo Garc\'ia-Portugu\'es$^{3}$}
\footnotetext[1]{Department of Quantitative Methods, CUNEF Universidad (Spain).}
\footnotetext[2]{Department of Mathematics ``Tullio Levi Civita'', University of Padova (Italy).}
\footnotetext[3]{Department of Statistics, Universidad Carlos III de Madrid (Spain).}
\footnotetext[4]{Corresponding author. e-mail: \href{mailto:abel.guada@cunef.edu}{abel.guada@cunef.edu}.}
\maketitle
%-----------------------------------------------%

\begin{abstract}
	We make a rigorous analysis of the existence and characterization of the free boundary related to the optimal stopping problem that maximizes the mean of an Ornstein--Uhlenbeck bridge. The result includes the Brownian bridge problem as a limit case. The methodology hereby presented relies on a time-space transformation that casts the original problem into a more tractable one with an infinite horizon and a Brownian motion underneath. We comment on two different numerical algorithms to compute the free-boundary equation and discuss illustrative cases that shed light on the boundary's shape. In particular, the free boundary generally does not share the monotonicity of the Brownian bridge case.
\end{abstract}
\begin{flushleft}
\end{flushleft}

%-------------------------------------------------%
\section{Introduction}\label{sec:intro}
%-------------------------------------------------%

Since their first appearance in the seminal monograph of \cite{Wald_1947_sequential}, Optimal Stopping Problems (OSPs) have become ubiquitous tools in mathematical finance, stochastic analysis, and mathematical statistics, among many other fields. Particularly, OSPs that are non-homogeneous in time are known to be mathematically challenging and, compared to the time-homogeneous counterpart, the literature addressing this topic is scarce, non-comprehensive, and often heavy on smoothness conditions. Markov bridges are not only time-inhomogeneous processes, but they also fail to meet the common assumption of Lipschitz continuity of the underlying drift (see, e.g., \citet[Chapter 3]{Krylov_1980_controlled}, or \cite{Jacka_1992_finite-horizon}), as their drifts explode when time approaches the horizon, thus inherently adding an extra layer of complexity.

The first result in OSPs with Markov bridges was given by \cite{Shepp_1969_explicit}, who circumvented the complexity of dealing with a Brownian Bridge (BB) by using a time-space transformation that allowed reformulating the problem into a more tractable one with a Brownian motion underneath. Since then, more than fifty years ago, the use of Markov bridges in the context of OSPs has been narrowed to extending the result of \cite{Shepp_1969_explicit}: \cite{Ekstrom_2009_optimal} and \cite{Ernst_2015_revisiting} studied alternative methods of solutions; \cite{Ekstrom_2009_optimal} and \cite{DeAngelis_2020_optimal} looked at a broader class of gain functions, \cite{Glover_2020_optimally} randomized the horizon while \cite{Follmer_1972_optimal}, \cite{Leung_2018_optimal}, and \cite{Ekstrom_2020_optimal} analyzed the randomization of the bridge's terminal point.

In finance, the use of a BB in OSPs has been motivated by several applications. \cite{Boyce_1970_stopping} applied it to the optimal selling of bonds; \cite{Baurdoux_2015_optimal} suggested the use of a BB to model mispriced assets that could rapidly return to their fair price, or perishable commodities that become useless after a given deadline; and \cite{Ekstrom_2009_optimal} used a BB to model the \textit{stock-pinning} effect, that is, the phenomenon in which the price of a stock tends to be pulled towards the strike price of one of its underlying options with massive trading volumes at the expiration date. While these motivations encourage the investor to rely on a model with added information at the horizon, none of them are exclusive to a BB, its usage being rather driven by tractability issues. Thus, in those same scenarios, other bridge processes could be more appealing than the over-simplistic BB. In particular, we drive our attention to an Ornstein--Uhlenbeck Bridge (OUB) process, since its version without a fixed terminal point, the Ornstein--Uhlenbeck (OU) process, is often the reference model in many financial problems.

Indeed, OU processes are a go-to in finance when it comes to modeling assets with prices that fluctuate around a given level. This mean-reverting phenomenon has been systematically observed in a wide variety of markets. A good reference for either theory, applications, or empirical evidence of mean-reverting problems is \cite{Leung_2015_optimal_book}. An example is given by the pair trading strategy, which consists of holding a position in one asset as well as the opposite position in another, both assets known to be correlated in a way that the spread between their prices shows mean reversion. Recently, many authors have tackled pair trading by using an OSP approach with an OU process. \cite{Ekstrom_2011_optimal} found the best time to liquidate the spread in the presence of a stop-loss level; \cite{Leung_2015_optimal_paper} used a discounted double OSP to compute the optimal buy-low-sell-high strategy in a perpetual frame; and \cite{Kitapbayev_2017_optimal} extended that result to a finite horizon and took the viewpoint of investors entering the spread either buying or shorting.   

In this paper we solve the finite-horizon OSP featuring the identity as the gain function and an OUB as the underlying process. The solution is provided in terms of a non-linear, Volterra-type integral equation. Similarly to \cite{Shepp_1969_explicit}, our methodology relies on a time-space change that casts the original problem into an infinite-horizon OSP with a Brownian motion as the underlying process. Due to the complexity of our resulting OSP, we use a direct approach to solve it rather than using the common candidate-verification scheme. We then show that one can either apply the inverse transformation to recover the solution of the original OSP or, equivalently, solve the Volterra integral equation reformulated back in terms of OUB. It is worth highlighting that the BB framework is included in our analysis as a limit case. 

The rest of the paper is structured as follows. Section \ref{sec:formulation} introduces the main problem and some useful notation. In Section \ref{sec:reformulation} we derive the transformed OSP and establish its equivalence to the original one. The most technical part of the paper is relegated to Section \ref{sec:sol_reformulated}, in which we derive the solution of the reformulated OSP. From it, we use the reverse transformation to get the solution back to the original OSP in Section \ref{sec:sol_original}, where we also remark that both a BB and an OUB with general pulling level and terminal time are immediate consequences of our results. An algorithm for numerical approximations of the solution is given in Section \ref{sec:numerical_results}, along with a compendium of illustrative cases for different values of the OUB's parameters. Concluding remarks are relegated to Section \ref{sec:conclusions}.

%-------------------------------------------------%
\section{Formulation of the problem}\label{sec:formulation}
%-------------------------------------------------%

Let $X = \{X_t\}_{t \in [0, 1]}$ be an OUB with terminal value $X_1 = z$, $z\in\R$, and defined in the filtered space $(\Omega, \mathcal{F}, \Pr, \{\cF_t\}_{t \in [0, 1]})$. That is, for an OU process $\wt{X} = \{\wt{X}_t\}_{t \in [0, 1]}$, take $X$ such that $\mathrm{Law}(X, \Pr) = \mathrm{Law}(\wt{X}, \wt{\Pr}_{z}),$ where $\wt{\Pr}_{z} := \mathbb{P}\big(\cdot | \wt{X}_1 = z\big)$. It is well known (see, e.g., \cite{Barczy_2013_sample}) that $X$ is the unique strong solution of the Stochastic Differential Equation (SDE)
\begin{align}\label{eq:OUB_SDE}
    \rmd X_t = \mu(t, X_t)\rmd t + \gamma\rmd B_t,\ \ 0 \leq t \leq 1,
\end{align}
with $\gamma > 0$ and
\begin{align}\label{eq:OUB_drift}
	\mu(t, x) = \alpha \frac{z - \cosh(\alpha(1 - t)) x}{\sinh(\alpha(1 - t))},\quad \alpha \neq 0.
\end{align}
Note that we can take $\{\mathcal{F}_t\}_{t \in [0, 1]}$ as the natural filtration of the underlying standard Brownian motion $\{B_s\}_{t \in [0, 1]}$ in \eqref{eq:OUB_SDE}. The parameter $\alpha$ regulates the strength with which the OU process $\wt{X}$ is attracted to ($\alpha > 0$) or repulsed from ($\alpha < 0$) the origin, resulting in the strength of attraction of the OUB $X$ to the curve $t \mapsto z/\cosh(\alpha(1 - t))$.
	
Consider the finite-horizon OSP
\begin{align}\label{eq:OSP_OUB}
	V(t, x) := \sup_{\tau \leq 1 - t}\Es{X_{t + \tau}}{t, x},
 \end{align}
where $V$ is the value function and $\E_{t, x}$ represents the expectation under the probability measure $\Pr_{t, x}$ defined as $\Pr_{t, x}(\cdot) := \Pr(\cdot | X_t = x)$. The supremum above is taken under all random times $\tau$ in the underlying filtration, such that $t + \tau$ is a stopping time in $\{\mathcal{F}_t\}_{t \in [0, 1]}$. Henceforth, we will call $\tau$ a stopping time while keeping in mind that $t + \tau$ is the actual stopping time.

%-------------------------------------------------%
\section{Reformulation of the problem}\label{sec:reformulation}
%-------------------------------------------------%

\cite{Barczy_2013_sample} provide the following space-time transformed representation for $X$:
\begin{align*}
    X_t &= a_1(t, X_0, z) + a_2(t)B_{\psi(t)},
\end{align*}
where the functions $a_1$ and $a_2$ take the form
\begin{align*}
    a_1(t, x, z) := x\frac{\sinh(\alpha(1 - t))}{\sinh(\alpha)} + z\frac{\sinh(\alpha t)}{\sinh(\alpha)}, \quad a_2(t) := \gamma e^{\alpha t}\frac{\kappa(1) - \kappa(t)}{\kappa(1)},
\end{align*}
and $\psi:[0, 1)\rightarrow \R_+$ is the time transformation
$
    \psi(t) := \kappa(t)\kappa(1)/(\kappa(1) - \kappa(t)),
$
with $\kappa(t) := (2\alpha)^{-1}(1 - e^{-2\alpha t})$. Notice that
$
    t = \kappa^{-1}\lrp{\psi(t)\kappa(1)/(\psi(t) + \kappa(1))},
$
where $\kappa^{-1}(s) = -(2\alpha)^{-1}\ln(1 - 2\alpha s)$. The following identities can be easily checked:
\begin{align*}
    a_1(t, x, z) = \lrp{x + z\frac{\psi(t)e^{-\alpha}}{\kappa(1)}}\frac{1}{f\lrp{\frac{\psi(t)e^{-\alpha}}{\kappa(1)}}}, \quad a_2(t) = \frac{\gamma}{f\lrp{\frac{\psi(t)e^{-\alpha}}{\kappa(1)}}},
\end{align*}
with
\begin{align}\label{eq:f}
    f(s) := \sqrt{\lrp{e^{\alpha} + s}\lrp{e^{-\alpha} + s}}. 
\end{align}
Therefore, if we set the time change $s = \psi(t)e^{-\alpha}/\kappa(1 - u)$, we get the space change
\begin{align}
    X_t = \frac{X_0 + zs}{f(s)} + \frac{\gamma}{f(s)}B_{s\kappa(1)e^{\alpha}}
    = \frac{zs + \gamma\sqrt{\kappa(1)e^{\alpha}}}{f(s)}\lrp{B_{s} + \frac{X_0}{\gamma\sqrt{\kappa(1)e^{\alpha}}}}. \label{eq:OUB_to_BM}
\end{align}

Let $Y = \lrb{Y_s}_{s\geq0}$ be a Brownian motion starting at $Y_0 = y$ under the probability measure $\Pr_{y}$, that is, $\Pro{Y_0 = y}{y} = 1$. Consider the infinite-horizon OSP
\begin{align}\label{eq:OSP_BM}
    W_{c}(s, y) := \sup_{\sigma}\Es{G_{c}(s + \sigma, Y_\sigma)}{y},
\end{align}
with gain function
\begin{align}\label{eq:gain}
    G_{c}(s, y) := \frac{cs + y}{f(s)}
\end{align}
and $c \in \R$. The parameter $c$ can be interpreted as a discount ($c < 0$) or interest ($c > 0$) rate per unit of time. The operator $\E_{y}$ emphasizes that we are taking the mean with respect to $\Pr_{y}$, and the supremum in \eqref{eq:OSP_BM} is taken over all the stopping times $\sigma$ in the natural filtration of $\lrb{Y_s}_{s\geq0}$.

Solving an OSP means giving a tractable expression for the value function and finding a stopping time in which the supremum is attained. Thereby, we show in the next proposition the equivalence between \eqref{eq:OSP_OUB} and \eqref{eq:OSP_BM}, by providing formulae that relate $V$ to $W$, and switch from a stopping time that is optimal in the former problem (if it exists) to one optimal in the latter. 

\begin{proposition}[Time-space equivalence]\label{pr:OSP_equiv}\ \\
    Consider the time change $\upsilon:[0,1]\rightarrow \R$ such that $\upsilon(t) := \psi(t)e^{-\alpha}/\kappa(1)$, and the space change $\eta:\R\rightarrow \R$ with $\eta(x) := x/(\gamma\sqrt{\kappa(1)e^{\alpha}})$. Take $(t, x)\in[0, 1)\times\R$ and set $s = \upsilon(t)$, $c_z = \eta(z)$, and $y = \eta(x)$. Then:
    \begin{enumerate}[label=(\roman{*}), ref=(\textit{\roman{*}})]
        \item \label{pr:OSP_value_equiv} The following equation holds:
        \begin{align}\label{eq:value_equiv}
            V(t, x) = \frac{z}{c_z}W_{c_z}\lrp{s, y}. 
        \end{align}
        \item \label{pr:OSP_OST_equiv} The stopping time $\sigma^*(s, y)$ is optimal in \eqref{eq:OSP_BM} under $\Pr_{y}$ for $c = c_z$ if and only if
        \begin{align}\label{eq:OST_transform}
            \tau^*(t, x) := \upsilon^{-1}\lrp{\sigma^*(s, y)}
        \end{align}
        is optimal in \eqref{eq:OSP_OUB} under $\Pr_{t, x}$.
    \end{enumerate}
\end{proposition}

\begin{proof}
    \ref{pr:OSP_value_equiv}\ We have already proved this part of the proposition. Indeed, \eqref{eq:value_equiv} follows trivially from \eqref{eq:OSP_OUB} and \eqref{eq:OUB_to_BM}--\eqref{eq:gain}.
    
    \ref{pr:OSP_OST_equiv}\ Suppose that $\sigma^* = \sigma^*(s, y)$ is optimal in \eqref{eq:OSP_BM} under $\Pr_{y}$ for $c = c_z$. Assume that there exists a stopping time $\tau' = \tau'(t, x)$ that outperforms $\tau^* = \tau^*(t, x)$ defined in \eqref{eq:OST_transform}, and set $\sigma' = \sigma'(s, y) := \upsilon^{-1}(\tau')$. Then, by relying on \eqref{eq:OUB_to_BM}, we get that
    \begin{align*}
        \Es{G_{c_z}\lrp{s + {\sigma'}, Y_{\sigma'}}}{y} = \Es{X_{t + \tau'}}{t, x} > \Es{X_{t + \tau^*}}{t, x} = \Es{G_{c_z}\lrp{s + {\sigma^*}, Y_{\sigma^*}}}{y},
    \end{align*}
    which contradicts the fact that $\sigma^*$ is optimal in \eqref{eq:OSP_BM}. Then, we have proved the \textit{only if} part of the statement. By following similar arguments, one can prove that if $\sigma^*$ is suboptimal, so it is $\tau^*$, which proves the \textit{if} direction. 
\end{proof}

%-------------------------------------------------%
\section{Solution of the reformulated problem: a direct approach}\label{sec:sol_reformulated}
%-------------------------------------------------%

In this section we work out a solution for the OSP \eqref{eq:OSP_BM}.
For the sake of briefness and since there is no risk of confusion,  throughout the section we will 
use the notations $W = W_c$ and $G = G_c$, so that \eqref{eq:OSP_BM} can be rewritten as
\begin{align}\label{eq:OSP_BM_re}
    W(s, y) = \sup_{\sigma}\Es{G(s + \sigma, Y_\sigma)}{y}.
\end{align}

Notice that $0 \leq s/f(s) \leq 1$ and $f(s) \geq \sqrt{1 + s^2}$ for all $s\in\R_+$, $f(0) = 1$, and $f$ is increasing. Hence, the following holds for $M := \Esp{\sup_{0\leq u\leq1}\left|B_u\right|}$ and all $(s, y) \in\R_+\times\R$:
\begin{align}
    \Es{\sup_{u\geq 0
    }\left|G\lrp{s + u, Y_u}\right|}{y} &\leq |c| + \Es{\sup_{u\geq 0}\frac{\left|Y_u\right|}{f(u)}}{y} \nonumber
    \leq |c| + |y| + \Esp{\sup_{u\geq0}\frac{\left|B_u\right|}{\sqrt{1 +  u^2}}} \nonumber \\
    &\leq |c| + |y| + M + \Esp{\sup_{u\geq1}\frac{\left|B_u\right|}{\sqrt{1 +  u^2}}} \nonumber \\
    &= |c| + |y| + M + \Esp{\sup_{u\geq1}\frac{u}{\sqrt{1 +  u^2}}\left|B_{1/u}\right|} \nonumber \\
    &\leq |c| + |y| + M + \Esp{\sup_{u\geq1} \left|B_{1/u}\right|}
    = |c| + |y| + 2M, \label{eq:sup_bound}
\end{align}
where we used the time-inversion property of a Brownian motion in the first equality. Thereby, since $M<\infty$ (see, e.g., Identity 1.1.3 from \cite{Borodin_2002_handbook}) and $G$ is continuous, we get that (see, e.g., Corollary 2.9, Remark 2.10, and equation (2.2.80) in \cite{Peskir_2006_optimal}) the first hitting time 
\begin{align}\label{eq:OST}
    \sigma^*(s, y) = \inf\lrb{u \geq 0 : (s + u, Y_u)\in \cD}
\end{align} 
into the stopping set $\cD := \lrb{W = G}$ is optimal for \eqref{eq:OSP_BM_re}. That is,
\begin{align}\label{eq:value}
    W(s, y) = \Es{G\lrp{s + \sigma^*(s, y), Y_{\sigma^*(s, y)}}}{y}.
\end{align}
After applying Itô's lemma to $G\lrp{s + u, Y_u^y}$, substituting $u$ for $\sigma$ and $\sigma^*(s, y)$ respectively, and taking expectation (which cancels the martingale term), we get from \eqref{eq:OSP_BM} and \eqref{eq:value} the following alternative representations of~$W$:
\begin{align}\label{eq:value_ito}
        W(s, y) - G(s, y) &= \sup_{\sigma}\Es{\int_0^\sigma \bbL G\lrp{s + u, Y_u}\,\rmd u}{y} = \Es{\int_0^{\sigma^*(s, y)} \bbL G\lrp{s + u, Y_u}\,\rmd u}{y},
\end{align}
where $\bbL = \partial_t + \frac{1}{2}\partial_{xx}$ is the infinitesimal generator of $\lrb{\lrp{s + u, Y_u}}_{u\geq 0}$. Here and thereafter, $\partial_t$ and $\partial_x$ will stand, respectively, for the differential operator with respect to time and space, while $\partial_{xx}$ is a shorthand for $\partial_x\partial_x$. Note that $\bbL G = \partial_t G$. Since many of the proofs rely on the first-order partial derivatives of the gain function, we display them next for quick reference:
\begin{align}
    \partial_t G(s, y) &= \frac{c\lrp{f(s) - sf'(s)} - f'(s)y }{f^2(s)}, \label{eq:G_t} \\
    \partial_x G(s, y) &= \frac{1}{f(s)}. \label{eq:G_x}
\end{align}

To keep track of the initial condition in a way that does not change the underlying probability measure, we introduce the process $Y^y = \lrb{Y_s^y}_{s\geq 0}$ such that 
\begin{align*}
    \mathrm{Law}\lrp{\lrb{Y_s^y}_{s\geq0}, \Pr} = \mathrm{Law}\lrp{\lrb{Y_s}_{s\geq 0}, \Pr_y}.	
\end{align*}

Notice that the characterization of the Optimal Stopping Time (OST) in \eqref{eq:OST} is too abstract to work with. In the next proposition we characterize $\sigma^*(s, y)$ by means of a function called the Optimal Stopping Boundary (OSB), which is the frontier between $\cD$ and its complement $\cC := \lrb{W > G}$. We also derive some properties about the shape of the OSB that shed light on the geometry of $\cD$ and~$\cC$.

\begin{proposition}[Existence and shape of the optimal stopping boundary]\label{pr:boundary_existence}\ \\
    There exists a function $b:\R_+\rightarrow\R$ such that $\cD = \lrb{(s, y) : y \geq b(s)}$. Moreover, $c(f(s) - sf'(s))/f(s) < b(s) < \infty$ for all $s\in\R_+$.
\end{proposition}

\begin{proof}
    The claimed shape for the stopping set, $\cD = \lrb{(s, y) : y \geq b(s)}$, is a straightforward consequence of the fact that $y \mapsto (W - G)(s, y)$ is decreasing for all $s\in\R_+$, which follows after \eqref{eq:f}, \eqref{eq:value_ito}, and \eqref{eq:G_t}.
    
    We now see that $b(s) > c(f(s) - sf'(s))/f(s)$ for all $s>0$. Fix a pair $(s, y)$ such that \mbox{$\partial_t G(s, y) > 0$}. Then, the continuity of $\partial_t G$ allows to pick a ball $\cB$ such that $(s, y)\in\cB$ and \mbox{$\partial_t G > 0$ in $\cB$}. After recalling \eqref{eq:value_ito} and setting $\sigma_\cB$ as the first exit time of $\lrb{\lrp{s + u, Y_u^y}}_{u\geq 0}$ from $\cB$, we get that
    \begin{align*}
    W(s, y) - G(s, y) \geq \Es{\int_0^{\sigma_\cB} \partial_t G\lrp{s + u, Y_u}\,\rmd u}{y} > 0.
    \end{align*}
    We conclude then that $(s, y)\in\cC$. Finally, the claimed lower bound for $b$ comes after using \eqref{eq:G_t} to realize that $\partial_t G(s, y) > 0$ if and only if $y < c(f(s) - sf'(s))/f(s)$.
    
    We now prove $b(s) < \infty$ for all $s>0$. Let $X =
    \big\{X_t\big\}_{t \in [0, 1]}$ and $\wt{X} =
    \big\{\wt{X}_t\big\}_{t \in [0, 1]}$ be an OUB and a BB, respectively, with pinning point $X_1 = \wt{X}_1 = z$. The drift of $\wt{X}$ has the form $\wt{\mu}(t, x) = (z - x)/(1 - t)$. Define $m_z:[0, 1)\rightarrow\R$ such that
    \begin{align*}
    m_z(t) = z\frac{\sinh(\alpha(1 - t)) - \alpha(1 - t)}{\sinh(\alpha(1 - t)) - \alpha(1 - t)\cosh(\alpha(1 - t))},
    \end{align*}
    and notice that $\mu(t, x) \leq \wt{\mu}(t, x)$ if and only if $x \geq m_z(t)$. Take $\overline{M}_z := \sup_{t\in[0, 1)}m_z(t) < \infty$ and notice the following relation:
    \begin{align}\label{eq:Xt-bound}
     X_t \leq m_z(t) + |X_t - m_z(t)| \leq m_z(t) + |\wt{X}_t - m_z(t)| \leq \overline{M}_z + |\wt{X}_t - \overline{M}_z|.
    \end{align}
    The second inequality in \eqref{eq:Xt-bound} holds since the drift of the process $t\mapsto m_z(t) + |X_t - m_z(t)|$ is lower than the drift of $t\mapsto m_z(t) + |\wt{X}_t - m_z(t)|$ and, therefore, we can ensure that, pathwise, the first process is lower than the last one $\Pr$-a.s. (see Corollary 3.1 by \cite{Peng_2006_necessary}). Indeed, for $\varepsilon > 0$, define the function
    \begin{align*}
        g_\varepsilon(x, m) := 
        \begin{cases}
        |x - m|, & \text{ if } |x - m| \geq \varepsilon, \\
        \frac{1}{2}\lrp{\varepsilon + \varepsilon^{-1}(x - m)^2}, & \text{ if } |x - m| < \varepsilon
        \end{cases}
    \end{align*}
    and the processes
    \begin{align*}
        Y_t^{(1), \varepsilon} := m_z(t) + g_\varepsilon\lrp{X_t, m_z(t)}, \quad Y_t^{(2),\varepsilon} := m_z(t) + g_\varepsilon\lrp{\wt{X}_t, m_z(t)}.
    \end{align*}
    Denote by $\mu^{(i),\varepsilon}$ the drift of $Y^{(i), \varepsilon}$, $i = 1, 2$. We obtain the following after a straightforward use of the Itô formula:
    \begin{align*}
        (\mu^{(1),\varepsilon} - \mu^{(2),\varepsilon})(t, x) = 
        \partial_1 g_\varepsilon(x, m(t))(\mu(t, x) - \wt{\mu}(t, x)).
    \end{align*}
    Therefore, after recalling the definition of $m_z(t)$ and noticing that $x \mapsto g_\varepsilon(x, m(t))$ decreases for $x \leq m(t)$ and increases for $x\geq m(t)$, a direct use of Corollary 3.1 from \cite{Peng_2006_necessary} yields that $Y_t^{(1), \varepsilon} \geq Y_t^{(2), \varepsilon}$ for all $t\in[0, T]$ $\Pr$-a.s., which implies the claimed result after taking $\varepsilon \rightarrow 0$ and realizing that, in such a case, $g_\varepsilon(c, x) \downarrow |c - x|$.
    
    The third inequality in \eqref{eq:Xt-bound} is straightforward from the definition of $\overline{M}_z$. 
    Therefore, if we consider the OSP
    \begin{align*}
    \wt{V}_{\overline{M}_z}(t, x) = \sup_{\tau \leq 1 - t}\Es{\overline{M}_z + |\wt{X}_{t + \tau} - \overline{M}_z|}{t, x},
    \end{align*}
    we are allowed to state that $V \leq \wt{V}_{\overline{M}_z}$. If we take a pair $(t, x)\in[0, 1]\times[\overline{M}_z, \infty)$ within the stopping set related to $\wt{V}_{\overline{M}_z}$, then $V(t, x) \leq \wt{V}_{\overline{M}_z} = x$, meaning that $(t, x)$ lies in the stopping set of $V$. 
    Since it is known that the OSB related to $\wt{V}_{\overline{M}_z}$ is finite (actually, this is one of the few cases in which the explicit form of the OSP with a finite horizon is available; see, e.g., Theorem 3.2 in \cite{Ekstrom_2009_optimal}), so is the one related to $V$. Then, using \eqref{eq:value_equiv}, we conclude that $b$ is bounded from above. 
\end{proof}

We next show that $W$ is Lipschitz continuous on sets of the type $\R_+\times \cR$, where $\cR$ stands for a compact set in $\R$.

\begin{proposition}[Local Lipschitz continuity of the value function]\label{pr:W_Lipschitz}\ \\
    For any compact set $\cR\subset\R$, there exists a constant $L_\cR > 0$ such that
    \begin{align*}
        \left|W(s_1, y_1) - W(s_2, y_2)\right| \leq L_{\cR}\lrp{|s_1 - s_2| + |y_1 + y_2|},
    \end{align*}
    for all $(s_1, y_1), (s_2, y_2) \in \R_+\times\cR$.
\end{proposition}

\begin{proof}
    Take $(s_1, y_1), (s_2, y_2) \in \R_+\times\cR$ and realize that
    \begin{align*}
        W(s_1, y_1) - W(s_2, y_2) =&\; \sup_{\sigma}\Es{G(s_1 + \sigma, Y_\sigma)}{y_1} - \sup_{\sigma}\Es{G(s_1 + \sigma, Y_\sigma)}{y_2} \\
        &+ \sup_{\sigma}\Es{G(s_1 + \sigma, Y_\sigma)}{y_2} - \sup_{\sigma}\Es{G(s_2 + \sigma, Y_\sigma)}{y_2}.
    \end{align*}
    Notice from \eqref{eq:G_t} that the following relation holds:
    \begin{align*}
        \left|\partial_t G(s, y)\right| \leq K \lrp{1 + \frac{|y|}{f(u)}}.
    \end{align*}
    Then, since $|\sup_\sigma a_\sigma - \sup_\sigma b_\sigma|\leq \sup_\sigma|a_\sigma - b_\sigma|$, alongside Jensen's inequality, and \eqref{eq:G_t} and \eqref{eq:G_x}, we get that
    \begin{align*}
        \Big|\sup_{\sigma}\Es{G(s_1 + \sigma, Y_\sigma)}{y_1}& - \sup_{\sigma}\Es{G(s_1 + \sigma, Y_\sigma)}{y_2}\Big| \\
         \leq&\; \sup_{\sigma}\Esp{\left|G(s_1 + \sigma, Y_\sigma^{y_1}) - G(s_1 + \sigma, Y_\sigma^{y_2})\right|} \\
        =&\; \sup_{\sigma}\Esp{\frac{\left|Y_\sigma^{y_1} - Y_\sigma^{y_2}\right|}{f(s_1 + \sigma)}} = \frac{|y_1 - y_2|}{f(s_1)} \leq |y_1 - y_2|,
    \end{align*}
    where, in the last equality, we used the fact that $f$ increases and the representation $Y_s^y = y + B_s$, for the standard Brownian motion $\lrc{B_s}_{s\in\R_+}$. Likewise
    \begin{align*}
        \Big|\sup_{\sigma}\Es{G(s_1 + \sigma, Y_\sigma)}{y_2}& - \sup_{\sigma}\Es{G(s_2 + \sigma, Y_\sigma)}{y_2}\Big| \\
        \leq&\; \sup_{\sigma}\Esp{\left|G(s_1 + \sigma, Y_\sigma^{y_2}) - G(s_2 + \sigma, Y_\sigma^{y_2})\right|} \\
        =&\; |s_1 - s_2|\sup_{\sigma}\Esp{\left|\partial_t G(\xi, Y_\sigma^{y_2})\right|} \\
        \leq&\; |s_1 - s_2|K\lrp{1 + \Esp{\sup_{s\geq0}\frac{\left| Y_s^{y_2}\right|}{f(s)}}},
    \end{align*}
    where $\xi\in\lrp{\min\lrb{s_1, s_2}, \max\lrb{s_1, s_2}}$ is a random variable that follows from the mean value theorem. Since we already proved in \eqref{eq:sup_bound} that $\Esp{\sup_{s\geq0}\left| Y_s^{y_2}\right|/f(s)} < \infty$, the Lipschitz continuity of $W$ in $\R_+\times\cR$ follows. 
\end{proof}

Beyond Lipschitz continuity, it turns out that the value function attains a higher smoothness away from the boundary. While this assertion is trivial in the interior of the stopping region, where $W = G$, we prove in the next proposition that it also holds in the continuation set. In addition, we show that $\bbL W$ vanishes in $\cC$, which establishes the equivalence between \eqref{eq:OSP_BM_re} and a free-boundary problem.

\begin{proposition}[Higher smoothness of the value function and the free-boundary problem]\label{pr:W_smoothness}\ \\
    $W\in C^{1, 2}(\cC)$ and $\bbL W = 0$ in $\cC$.
\end{proposition}

\begin{proof}
    The fact that $\mathbb{L}W = 0$ in $\cC$ comes right after the strong Markov property of $\lrb{\lrp{s + u, Y_u}}_{u\geq0}$; see \citet[Section 7.1]{Peskir_2006_optimal} for more details.
	
	Since $W$ is continuous on $\cC$ (see Proposition \ref{pr:W_Lipschitz}) and the coefficients in the parabolic operator $\bbL = \partial_t + \frac{1}{2}\partial_{xx}$ are constant (it suffices to require local $\alpha$-H\"older continuity), then standard theory from parabolic partial differential equations \citep[Section 3, Theorem 9]{Friedman_1964_partial} guarantees that, for an open rectangle $R\subset \cC$, the first initial-boundary value problem
    \begin{subequations} \label{eq:PDE}
    \begin{align}
	\bbL g &= 0 &&\hspace*{-3.5cm} \text{in } R, \label{eq:PDE1} \\
        g &= W &&\hspace*{-3.5cm} \text{on } \partial R \label{eq:PDE2}
    \end{align}
    \end{subequations}
    has a unique solution $g\in C^{1, 2}(R)$. Therefore, if we denote by $\tau_{R^c}$ the first time $(s + u, Y_u)$ exits $R$, %we can use Dynkin's formula on $g(s + u, Y_u)$ at $u = \tau_{R^c}$, that is, the first time $(s + u, Y_u)$ exits $R$, %and then take $\Pr_y$-expectation with $y \in R$, which guarantees the vanishing of the martingale term and yields, together with \eqref{eq:PDE1} and \eqref{eq:PDE2}, the equality 
    we obtain the following for $y \in R$:
    \begin{align*}
    W(s, y) = \Es{W(s + \tau_{R^c}, Y_{\tau_{R^c}})}{y} = \Es{g(s + \tau_{R^c}, Y_{\tau_{R^c}})}{y} = g(t, x) + \Es{\int_0^{\tau_{R^c}}\bbL g(s+u, Y_u)\,\rmd u}{y} = g(t, x),
    \end{align*}
    where the first equality is due to the strong Markov property, the second one relies on \eqref{eq:PDE2}, the third one comes after a straightforward application of Dynkin's formula, and \eqref{eq:PDE1} was used to obtain the last one.  %Finally, due to the strong Markov property, $\E_y[W(s + \tau_{R^c}, Y_{\tau_{R^c}})] = W(s, y)$. 
\end{proof}

Not only the gain function has continuous partial derivatives away from the boundary, but we can provide relatively explicit forms for those derivatives, as shown in the next proposition.

\begin{proposition}[Partial derivatives of the value function]\label{pr:W_t_&_W_x}\ \\
Let $\sigma^* = \sigma^*(s, y)$, for $(s, y)\in\cC$, and $a := e^{-\alpha} + e^{\alpha}$. Then, 
\begin{align}\label{eq:W_t}
        &\partial_t W(s, y) = \partial_t G(s, y) + \Esp{\int_s^{s + \sigma^*}\frac{1}{f^3(u)}\lrp{-c\lrp{a + 3u} +  \frac{3\lrp{a + 2u}^2}{4f^2(u)} - Y_{u - s}^y}\, \rmd u}
    \end{align}
    and
    \begin{align}\label{eq:W_x}
        \partial_x W(s, y) = \Esp{\frac{1}{f(s + \sigma^*)}}.
    \end{align}
\end{proposition}

\begin{proof}
    Take $(s, y)\in\cC$ and $\varepsilon > 0$. Due to \eqref{eq:OSP_BM_re} and \eqref{eq:value}, we get the following for $\sigma^* = \sigma^*(s, y)$:
    \begin{align*}
        \varepsilon^{-1}\lrp{W(s, y) - W(s - \varepsilon, y)} 
        &\leq \varepsilon^{-1}\Esp{G(s + \sigma^*, Y_{\sigma^*}^y) - G(s - \varepsilon + \sigma^*, Y_{\sigma^*}^y)}.
    \end{align*}
    Hence, by letting $\varepsilon\rightarrow 0$, using the dominated convergence theorem (see \eqref{eq:sup_bound}), and recalling that $W\in C^{1,2}(\cC)$ (see Proposition \ref{pr:W_smoothness}), we get that
    \begin{align}\label{eq:W_t<}
        \partial_t W(s, y) &\leq \Esp{\partial_t G(s + \sigma^*, Y_{\sigma^*}^y)} = \partial_t G(s, y) + \Esp{\int_s^{s + \sigma^*}\bbL\partial_t G(u, Y_{s - u}^y)\, \rmd u}.
    \end{align}
    In the same fashion we obtain 
    \begin{align*}
        \varepsilon^{-1}\lrp{W(s + \varepsilon, y) - W(s, y)} 
        &\geq \varepsilon^{-1}\Esp{G(s + \varepsilon + \sigma^*, Y_{\sigma^*}^y) - G(s + \sigma^*, Y_{\sigma^*}^y)}.
    \end{align*}
    Thus, by arguing as in \eqref{eq:W_t<} we get the reverse inequality, and therefore \eqref{eq:W_t} is proved after computing $\bbL\partial_t G(u, Y_{s - u}^y) = \partial_{tt} G(u, Y_{s - u}^y)$.
    
    To get the analog result for the space coordinate, notice that
    \begin{align*}
        \varepsilon^{-1}\lrp{W(s, y) - W(s, y - \varepsilon)} &\leq \varepsilon^{-1}\Esp{W(s + \sigma^*, Y_{\sigma^*}^y) - W(s + \sigma^*, Y_{\sigma^*}^{y - \varepsilon})} \\
        &\leq \varepsilon^{-1}\Esp{G(s + \sigma^*, Y_{\sigma^*}^y) - G(s + \sigma^*, Y_{\sigma^*}^{y - \varepsilon})} \\
        &= \Esp{\frac{1}{f(s + \sigma^*)}},
    \end{align*}
    while the same reasoning yields the inequality $\varepsilon^{-1}\lrp{W(s, y + \varepsilon) - W(s, y)} \geq \Esp{1/f(s + \sigma^*)}$, and then, by letting $\varepsilon\rightarrow 0$, we get \eqref{eq:W_x}. 
\end{proof}

So far we have proved that the solution of \eqref{eq:OSP_BM_re} also solves the free-boundary problem
\begin{subequations} \label{eq:free-boundary}
	\begin{align}
		\bbL W(s, y) &= 0 &&\hspace*{-3cm} \text{for } y < b(t) , \label{eq:free-boundary1}\\
		W(s, y) &> G(s, y) &&\hspace*{-3cm} \text{for } y < b(t), \label{eq:free-boundary2}\\
		W(s, y) &= G(s, y) &&\hspace*{-3cm} \text{for } 
		y\geq b(t). \label{eq:free-boundary3}
	\end{align}
\end{subequations}
However, an additional condition for the value function on the free boundary is required to guarantee a unique solution. Roughly speaking, this condition comes in the form of smoothly binding the value and the gain functions with respect to the space coordinate, provided that the optimal boundary is (probabilistically) regular for the interior of $D$, that is, if after starting at a point $(s, y) \in \partial\cC$, the process enters the interior of $D$ immediately $\Pr_y$-a.s. (see \cite{DeAngelis_2020_global}). This type of regularity can be derived for locally Lipschitz continuous OSBs (see Proposition \ref{pr:smooth-fit} ahead).

In the next proposition we show that the boundary is Lipschitz continuous on any bounded interval. The proof is inspired by Theorem 4.3 from \cite{DeAngelis_2019_Lipschitz}, which states the boundary's local Lipschitz continuity for high-dimensional processes with some regularity conditions. Our settings do not satisfy Assumption (D) in \cite{DeAngelis_2019_Lipschitz}, which establishes a relation between the partial derivatives of $G$.

\begin{proposition}[Lipschitz continuity of the optimal stopping boundary]\label{pr:Lipschitz_boundary}\ \\
For any closed interval $I := [\underline{s}, \overline{s}]\subset\R_+$, there exists a constant $L_I > 0$ such that
\begin{align}\label{eq:|b'|<}
    |b(s_1) - b(s_2)| \leq L_I,
\end{align}
whenever $s_1, s_2\in I$.
\end{proposition}
 
\begin{proof}
    Consider the function $H:I\times\R\rightarrow\R_+$, for a closed interval $I\subset\R_+$, defined as $H(s, y) = W(s, y) - G(s, y)$. Proposition \ref{pr:boundary_existence} entails that $b$ is bounded from below, and thus we can choose a constant $r\in \R$ such that $r < \inf\lrb{b(s) : s\in I}$. Since $I\times\lrb{r}\subset\cC$, $H$ is continuous (see Proposition \ref{pr:W_Lipschitz}) and $H|_{I\times\lrb{r}} > 0$. Then, there exists $a > 0$ such that $H(s, r) > a$ for all $s\in I$. Therefore, for all $\delta$ such that $0 < \delta \leq a$, the equation $H(s, y) = \delta$ has a solution in $\cC$ for all $s\in I$. Moreover, this solution is unique for each $s$ since $\partial_x H < 0$ in $\cC$ (see Proposition \ref{pr:W_t_&_W_x}), and we denote it by $b_\delta(s)$, where $b_\delta:I\rightarrow \R$. 
    Away from the boundary, $H$ is regular enough to apply the implicit function theorem that guarantees that $b_\delta$ is differentiable and
    \begin{align}\label{eq:b_delta'}
        b_\delta'(s) = -\partial_t H(s, b_\delta(s)) / \partial_x H(s, b_\delta(s)).
    \end{align}
    Note that $b_\delta$ is decreasing in $\delta$ and therefore converges pointwise to some limit function $b_0$, which satisfies $b_0 \leq b$ in $I$ as $b_\delta < b$ for all $\delta$. Since $H(s, b_\delta(s)) = \delta$ and $H$ is continuous, it follows that $H(s, b_0(s)) = 0$ after taking $\delta\rightarrow 0$, which means that $b_0 \geq b$ in $I$ and hence $b_0 = b$ in $I$.
	
    Take $(s, y)\in\cC$ such that $y>r$. Set $\sigma^* = \sigma^*(s, y)$ and consider
    \begin{align*}
        \sigma_r = \sigma_r(s, y) := \inf\lrb{u\geq 0 : \lrp{s + u, Y_u^y} \notin I\times(r, \infty)}.
    \end{align*} 
    Recalling \eqref{eq:W_t}, it is easy to check that %there exists a constant $K_I^{(1)} > 0$ such that
    \begin{align}\label{eq:H_t<_1}
        \left|\partial_t H(s, y)\right| \leq K_I^{(1)} m(s, y),
    \end{align}
    with
    \begin{align*}
        m(s, y) := \Es{\int_0^{\sigma^*}\lrp{1 + \frac{\left|Y_u\right|}{f^2(s + u)}}\,\rmd u}{y}
    \end{align*}
    and
    \begin{align*}
        K_I^{(1)} := \max\lrb{\sup_{u\in\R_+}\frac{1}{f(u)}, \sup_{u\in\R_+}\left|\frac{3\lrp{a + 2u}^2}{4f^5(u)} - \frac{c\lrp{a + 3u}}{f^3(u)}\right|}.        
    \end{align*}
    %Indeed, $K_I^{(1)} = \max\lrb{\sup_{u\in\R_+}\frac{1}{f(u)}, \sup_{u\in\R_+}\left|\frac{3\lrp{a + 2u}^2}{4f^5(u)} - \frac{c\lrp{a + 3u}}{f^3(u)}\right|}$, for example, verifies \eqref{eq:H_t<_1}.
    Using the tower property of conditional expectation, alongside the strong Markov property, we get
    \begin{align}
        m&(s, y) \nonumber\\
        &= \Es{\int_0^{\sigma^*\wedge\sigma_r}\lrp{1 + \frac{\left|Y_u\right|}{f^2(s + u)}}\,\rmd u + \mathbbm{1}\lrp{\sigma_r \leq \sigma^*}\int_{\sigma_r}^{\sigma^*}\lrp{1 + \frac{\left|Y_u\right|}{f^2(s + u)}}\,\rmd u}{y} \nonumber \\
        &= \Es{\int_0^{\sigma^*\wedge\sigma_r}\lrp{1 + \frac{\left|Y_u\right|}{f^2(s + u)}}\,\rmd u + \mathbbm{1}\lrp{\sigma_r \leq \sigma^*} \Es{\int_{\sigma_r}^{\sigma_r + \sigma^*\lrp{\sigma_r, Y_{\sigma_r}}}\lrp{1 + \frac{\left|Y_u\right|}{f^2(s + u)}}\,\rmd u \Big| \cF_{\sigma_r}}{y}}{y} \nonumber \\ 
        &= \Es{\int_0^{\sigma^*\wedge\sigma_r}\lrp{1 + \frac{\left|Y_u\right|}{f^2(s + u)}}\,\rmd u + \mathbbm{1}\lrp{\sigma_r \leq \sigma^*} \Es{\int_{0}^{\sigma^*\lrp{\sigma_r, Y_{\sigma_r}}}\lrp{1 + \frac{\left|Y_u\right|}{f^2(s + \sigma_r + u)}}\,\rmd u}{Y_{\sigma_r}}}{y} \nonumber \\
        &= \Es{\int_0^{\sigma^*\wedge\sigma_r}\lrp{1 + \frac{\left|Y_u\right|}{f^2(s + u)}}\,\rmd u + \mathbbm{1}\lrp{\sigma_r \leq \sigma^*} m(s + \sigma_r, Y_{\sigma_r})}{y}. \label{eq:m=}
    \end{align}
    Notice that, for $c < r < y < b(s)$, $(s + \sigma_r, Y_{\sigma_r}^y) \in \Gamma_s$ on the set $\lrb{\sigma_r \leq \sigma^*}$, with $\Gamma_s := \lrb{(s, \bar{s})\times\lrb{r}} \cup \lrb{\bar{s}\times[r, b(\bar{s}))}$ and $\bar{s}:= \sup\lrb{s: s \in I}$. Hence, the following holds true on the set $\lrb{\sigma_r \leq \sigma^*}$:
    \begin{align}
        m\lrp{s + \sigma_r, Y_{\sigma_r}^y} &\leq \sup_{(t, x) \in \Gamma_s}m\left(t, x\right) \nonumber \\
        &\leq \sup_{(t, x) \in \Gamma_s} \Es{\int_{0}^{\infty}\lrp{1 + \frac{\left|Y_u\right|}{f^2(t + u)}}\,\rmd u}{x} \nonumber \\
        &\leq \sup_{(t, x) \in \Gamma_s} \int_{0}^{\infty}\lrp{1 + \frac{\left|x\right|}{f^2(t + u)}}\,\rmd u +  \int_{0}^{\infty}\frac{\Esp{\left|B_u\right|}}{f^2(t + u)}\,\rmd u \nonumber \\
        &\leq \int_{0}^{\infty}\lrp{1 +\frac{ \left|b(\bar{s})\right|}{f^2(u)}}\,\rmd u + \int_{0}^{\infty}\sqrt{\frac{2}{\pi}}\frac{\sqrt{u}}{f^2(u)}\,\rmd u < \infty.  \label{eq:m_bound}
    \end{align}
    By plugging \eqref{eq:m_bound} into \eqref{eq:m=}, after observing that $\lrp{1 + \left|Y_u\right|/f^2(s + u)} \leq 1 + \max\lrb{|\sup_{s\in I}b(s)|, |r|}$, and recalling \eqref{eq:H_t<_1}, we obtain the following for some constant $K_I^{(2)} > 0$:
    \begin{align}\label{eq:H_t<_2}
        \left|\partial_t H(s, y)\right| \leq K_I^{(2)} \Es{\sigma_\delta\wedge\sigma_r + \mathbbm{1}\lrp{\sigma_r \leq \sigma_\delta}}{y}.
    \end{align}
    Arguing as in \eqref{eq:m=} and recalling \eqref{eq:G_x} along with \eqref{eq:W_x}, we get that
    \begin{align}
        \left|\partial_x H(s, y)\right| &= \Es{\frac{1}{f(s)} - \frac{1}{f(s + \sigma^*)}}{y} = \Es{\int_0^{\sigma^*}-\partial_t(1/f)(s + u)\,\rmd u}{y} \nonumber \\
        &= \Es{\int_0^{\sigma^*\wedge\sigma_r}-\partial_t(1/f)(s + u)\,\rmd u + \mathbbm{1}\lrp{\sigma_r \leq \sigma^*} \left|\partial_x H(s + \sigma_r, Y_{\sigma_r})\right|}{y} \nonumber \\
        &\geq \Es{\int_0^{\sigma^*\wedge\sigma_r}-\partial_t(1/f)(s + u)\,\rmd u + \mathbbm{1}\lrp{\sigma_r \leq \sigma^*, \sigma_r < \overline{s} - s}\left|\partial_x H(s + \sigma_r, r)\right|}{y}. \label{eq:H_x>_1}
    \end{align}
    Take $\varepsilon > 0$ such that $\cR_\varepsilon := [\underline{s}, \overline{s} + \varepsilon]\times\lrp{r-\varepsilon, r + \varepsilon}\subset \cC$, and consider the stopping time $\sigma_{\varepsilon} = \inf\lrb{u\geq0: Y_u^r \notin \cR_\varepsilon}$. Observe that $\sigma^*(s, r) > \sigma_{\varepsilon}$ for all $s\in I$. Then,
    \begin{align}
        \left|\partial_x H(s + \sigma_r, r)\right| &\geq \inf_{s\in I} \left|\partial_x H(s, r)\right| = \inf_{s\in I} \Es{\frac{1}{f(s)} - \frac{1}{f(s + \sigma^*(s, r))}}{r} \nonumber \\
        &\geq \inf_{s\in I} \Es{\frac{1}{f(s)} - \frac{1}{f(s + \sigma_\varepsilon)}}{r} \nonumber \\
        &\geq \inf_{s\in I} \lrp{\frac{1}{f(s)} - \frac{1}{f(\overline{s} + \varepsilon)}}\Pro{\sigma_\varepsilon = \overline{s} + \varepsilon - s}{r} \nonumber \\
        &= \lrp{\frac{1}{f(\overline{s})} - \frac{1}{f(\overline{s} + \varepsilon)}}\Pro{\sigma_\varepsilon = \overline{s} + \varepsilon - \underline{s}}{r} \nonumber \\
        &= \lrp{\frac{1}{f(\overline{s})} - \frac{1}{f(\overline{s} + \varepsilon)}}\Prob{\sup_{u\leq \overline{s} + \varepsilon - \underline{s}} |B_u| < \varepsilon} \label{eq:H_x>_2} > 0,
    \end{align}
    where we used the fact that $s\mapsto 1/f(s) - 1/f(s + u)$ is decreasing for all $u\geq 0$. The third inequality comes after using the law of the total probability and discarding the positive addend relative to the event $\lrb{\sigma_\varepsilon < \overline{s} + \varepsilon - s}$, while the last identity, coming after the equivalence of the two probability terms, is a direct consequence of the definition of $\sigma_\varepsilon$.
    After noticing that $-\partial_t(1/f)$ is positive and decreasing, which means that $-\partial_t(1/f)(s + u)\geq -\partial_t(1/f)(\overline{s}) > 0$ \mbox{for all $u\leq \sigma_r$}, and by plugging \eqref{eq:H_x>_2} into \eqref{eq:H_x>_1}, we obtain, for a constant $K_{I, \varepsilon}^{(3)} > 0$,
    \begin{align}
        \left|\partial_x H(s, y)\right| \geq K_{I}^{(3)}\Es{\sigma^*\wedge\sigma_r + \mathbbm{1}\lrp{\sigma_r \leq \sigma^*, \sigma_r < \overline{s} - s}}{y}. \label{eq:H_x>_3}
    \end{align}
    Therefore, using \eqref{eq:H_t<_2} and \eqref{eq:H_x>_3} in \eqref{eq:b_delta'} yields the following bound for some constant $K_{I}^{(4)} > 0$, $y_\delta = b_\delta(s)$, and $\sigma_\delta = \sigma^*(s, y_\delta)$:
    \begin{align}
        \left|b_\delta'(s)\right| &\leq  K_{I}^{(4)}\frac{\Es{\sigma_\delta\wedge\sigma_r + \mathbbm{1}\lrp{\sigma_r \leq \sigma_\delta}}{y_\delta}}{\Es{\sigma_\delta\wedge\sigma_r + \mathbbm{1}\lrp{\sigma_r \leq \sigma_\delta, \sigma_r < \overline{s} - s}}{y_\delta}} \nonumber \\
        &\leq K_{I}^{(4)}\lrp{1 + \frac{\Pro{\sigma_r \leq \sigma_\delta}{y_\delta}}{\Es{\sigma_\delta\wedge\sigma_r + \mathbbm{1}\lrp{\sigma_r \leq \sigma_\delta, \sigma_r < \overline{s} - s}}{y_\delta}}} \nonumber \\
        &\leq K_{I}^{(4)}\lrp{1 + \frac{\Pro{\sigma_r \leq \sigma_\delta, \sigma_r = \bar{s} - s}{y_\delta}}{\Es{\sigma_\delta\wedge\sigma_r}{y_\delta}} + \frac{\Pro{\sigma_r \leq \sigma_\delta, \sigma_r < \bar{s} - s}{y_\delta}}{\Es{\sigma_\delta\wedge\sigma_r + \mathbbm{1}\lrp{\sigma_r \leq \sigma_\delta, \sigma_r < \overline{s} - s}}{y_\delta}}} \nonumber \\
        &\leq K_I^{(4)}\lrp{2 + \frac{\Pro{\sigma_r \leq \sigma_\delta, \sigma_r = \bar{s} - s}{y_\delta}}{\Es{\mathbbm{1}\lrp{\sigma_r\leq\sigma_\delta, \sigma_r = \bar{s} - s}\lrp{\sigma_\delta\wedge\sigma_r}}{y_\delta}}} \nonumber \\
        &\leq K_I^{(4)}\lrp{2 + \frac{1}{\bar{s} - s}}. \label{eq:|b_delta'|<}
    \end{align}

    If we set $I_\varepsilon = [\underline{s}, \bar{s}  - \varepsilon]$ for $\varepsilon > 0$ small enough, then, by relying on \eqref{eq:|b_delta'|<}, we obtain the existence of a constant $L_{I_\varepsilon} > 0$, independent from $\delta$, such that $|b_\delta'(s)| < L_{I_\varepsilon}$ for all $s\in I_\varepsilon$ and $0 < \delta \leq a$. We are thus able to use the Arzelà--Ascoly's theorem to guarantee that $b_\delta$ converges to $b$ uniformly with respect to $\delta$ in $I_\varepsilon$. The proof is concluded after realizing that the Lipschitz continuity property is closed under the uniform limit operation, and that $\varepsilon$ can be chosen arbitrarily.  %Since $\varepsilon > 0$ and $I$ were arbitrarily chosen, we then conclude that $b$ is anywhere differentiable and \eqref{eq:|b'|<} holds true. 
\end{proof}

Once we have the Lipschitz continuity of the boundary on bounded sets, %this implying piecewise monotonicity, 
we proceed to illustrate in the following proposition how to obtain the principle of smooth fit, which, as we highlighted before, is required to provide a unique solution to the associated free-boundary problem \eqref{eq:free-boundary1}--\eqref{eq:free-boundary3}.

\begin{proposition}[The smooth-fit condition]\label{pr:smooth-fit}\ \\
    For all $s\geq 0$, $y\mapsto W(s, y)$ is differentiable at $y = b(s)$. Moreover, $\partial_x W(s, b(s)) = \partial_x G(s, b(s))$.
\end{proposition}

\begin{proof}
    Recall that we have already obtained in \eqref{eq:W_x} an explicit form for $\partial_x W$ away from the boundary, namely,
    \begin{align*}
    \partial_x W(s, y) = \Esp{\frac{1}{f(s + \sigma^*(s, y))}} \ , \quad (s,y) \in \cC. 
    \end{align*}
    The principle of smooth fit is just the validation of this formula whenever $y = b(s)$, $s\in\R_+$.
    
    We have that $\partial_x W(s, b(s)^+) = \partial_x G(s, b(s)) = 1/f(s)$, as $\sigma^*(s, y) = 0$ for all $y \geq b(s)$. By relying on the law of the iterated logarithm alongside the local Lipschitz continuity of $b$, we get that $(s, b(s))$ is probabilistically regular for the interior of $\cD$, that is,
    \begin{align*}
        %&\Prob{\inf\lrc{u > 0 : Y_u^{s, b(s)} < b(s+u)} = 0} \\
        %&= \lim_{\varepsilon\downarrow 0} \Prob{\inf\lrc{u > 0 : Y_u^{s, b(s)} < b(s+u)} < \varepsilon} 
        \lim_{\varepsilon\downarrow 0}\Prob{\inf_{u\in(0,\varepsilon)} \lrp{Y_u^{s, b(s)} - b(s+u)} < 0} &= \lim_{\varepsilon\downarrow 0}\Prob{\inf_{u\in(0,\varepsilon)} \frac{Y_u^{s, b(s)} - b(s+u)}{\sqrt{2u\ln(\ln(1/u))}} < 0} \\
        &\geq \lim_{\varepsilon\downarrow 0}\Prob{\inf_{u\in(0,\varepsilon)} \frac{Y_u^{s, b(s)} - b(s) + L_s u}{\sqrt{2u\ln(\ln(1/u))}} < 0}\\
        &= \Prob{\liminf_{u\downarrow 0} \frac{Y_u^{s, b(s)} - b(s) + L_s u}{\sqrt{2u\ln(\ln(1/u))}} < 0} = 1,
    \end{align*}
    for some $L_s > 0$. Corollary 6 from \cite{DeAngelis_2020_global} then provides
    %\citet[Corollary 8]{Cox_2015_embedding}, alongside the fact that our OSB is piecewise monotonic and continuous, 
    that $\sigma^*(s, b(s)^-) = \sigma^*(s, b(s)) = 0\ \Pr$-a.s. and, hence, the dominated convergence theorem entails that $\partial_x W(s, b(s)^-) = 1/f(s) = \partial_x G(s, b(s))$, thus concluding that the smooth-fit condition holds. 
\end{proof}

We are now in the position of getting a tractable characterization of both the value function and the OSB. Propositions \ref{pr:boundary_existence}--\ref{pr:smooth-fit} allow us to use an extension of Itô's lemma on the function $W(s + t, Y_t)$ for $t\geq 0$. This extension was originally derived by \cite{Peskir_2005_change} and later restated, in a way applies more directly to our framework, in Lemma A2 from \cite{DAuria_2020_discounted}. Recalling that $\bbL W = 0$ on $\cC$ and $W = G$ on $\cD$, and after taking $\Pr_y$-expectation (which cancels the martingale term), we get
\begin{align}
    W(s, y) &= \Es{W(s + t, Y_t)}{y} - \Es{\int_0^t (\bbL W)\lrp{s + u, Y_u}\,\rmd u}{y} \nonumber \\
    &= \Es{W(s + t, Y_t)}{y} - \Es{\int_0^t \partial_t G\lrp{s + u, Y_u}\mathbbm{1}\lrp{Y_u \geq b(s + u)}\,\rmd u}{y}, \label{eq:pricing_formula_aux}
\end{align}
where the local-time term does not appear due to the smooth-fit condition.

\begin{lemma}\label{lm:W->c}
    For all $(s, y)\in\R_+\times\R$,
    \begin{align*}
    \lim_{u\rightarrow\infty}\Es{W(s + u, Y_u)}{y} = c.
    \end{align*}
\end{lemma}

\begin{proof}
    The Markov property of $Y$, together with the fact that both $s\mapsto s/f(s)$ and $s\mapsto f(s)$ are increasing and $s/f(s)\rightarrow 1$ as $s\rightarrow\infty$, implies that  
    \begin{align}
        \Es{W(s + u, Y_u)}{y} &= \Es{\sup_{\sigma}\Es{G\lrp{s + u + \sigma, Y_\sigma}}{Y_u}}{y} \leq \Es{\Es{\sup_{r\geq 0}G\lrp{s + u + r, Y_r}}{Y_u}}{y} \nonumber \\
        &= \Es{\Es{\sup_{r\geq 0}\lrb{c\frac{s + u + r}{f(s + u + r)} + \frac{Y_r}{f(s + u + r)}}}{Y_u}}{y} \nonumber \\
        &\leq c\lrp{\mathbbm{1}(c > 0) + \frac{s + u}{f(s + u)}\mathbbm{1}(c \leq 0)} + \Es{\sup_{r\geq 0}\frac{Y_{u+ r}}{f(u + r)}}{y}, \label{eq:W->c_upper_bound}
    \end{align}
    and
    \begin{align}
        \Es{W(s + u, Y_u)}{y} &\geq
        \Es{\Es{\inf_{r\geq 0}G\lrp{s + u + r, Y_r}}{Y_u}}{y} \nonumber \\
        &\geq c\lrp{\mathbbm{1}(c < 0) + \frac{s + u}{f(s + u)}\mathbbm{1}(c \geq 0)} + \Es{\inf_{r\geq 0}\frac{Y_{u+ r}}{f(s + u + r)}}{y}. \label{eq:W->c_lower_bound}
    \end{align}
    Notice that
    \begin{align*}
        \lim_{u\rightarrow\infty}\Es{\sup_{r\geq 0}\frac{Y_{u+ r}}{f(u + r)}}{y} = \Es{\lim_{u\rightarrow\infty}\sup_{r\geq u}\frac{Y_{r}}{f(r)}}{y} = \Es{\limsup_{u\rightarrow\infty}\frac{Y_{u}}{f(u)}}{y} = 0,
    \end{align*}
    where in the first equality we applied the monotone convergence theorem and in the second one we used the law of the iterated logarithm as an estimate of the convergence of the process in the numerator. A similar argument yields 
    \begin{align*}
        \lim_{u\rightarrow\infty}\Es{\inf_{r\geq 0}\frac{Y_{u+ r}}{f(s + u + r)}}{y} = 
        \Es{\lim_{u\rightarrow\infty}\inf_{r\geq u}\frac{Y_{r}}{f(s + r)}}{y} = \Es{\liminf_{u\rightarrow\infty}\frac{Y_{u}}{f(s + u)}}{y} = 0. 
    \end{align*}
    Thus, we can take $u\rightarrow\infty$ in both \eqref{eq:W->c_upper_bound} and \eqref{eq:W->c_lower_bound} to complete the proof. 
\end{proof}

By taking $t\rightarrow\infty$ in \eqref{eq:pricing_formula_aux} and relying on Proposition \ref{lm:W->c}, we get the following pricing formula for the value function:
\begin{align}
    W(s, y) &= c - \Es{\int_0^\infty (\bbL W)\lrp{s + u, Y_u}\,\rmd u}{y} \nonumber \\
    &= c - \Es{\int_0^\infty \partial_t G\lrp{s + u, Y_u}\mathbbm{1}\lrp{Y_u \geq b(s + u)}\,\rmd u}{y}. \label{eq:pricing_formula}
\end{align}
We can obtain a more tractable version of \eqref{eq:pricing_formula} by exploiting the linearity of $y\mapsto \partial_t G(s, y)$ (see \eqref{eq:G_t}) as well as the Gaussianity of $Y_u$. Specifically, since $Y_u\sim \cN(y, u)$ under $\Pr_y$, then $\Es{Y_u\mathbbm{1}\lrp{Y_u\geq x}}{y} = \bar{\Phi}((x - y)/\sqrt{u})y + \sqrt{u}\phi((x - y)/\sqrt{u})$, where $\bar{\Phi}$ and $\phi$ denote the survival and the density functions of a standard normal random variable, respectively. By shifting the integrating variable $s$ units to the right, we get that
\begin{align}\label{eq:pricing_formula_refined}
    W(s, y) &= c - \int_s^\infty \frac{1}{f(u)}\lrp{c\bar{\Phi}_{s, y, u, b(u)} - \frac{(a + 2u)\lrp{(y + cu)\bar{\Phi}_{s, y, u, b(u)} + \sqrt{u - s}\phi_{s, y, u, b(u)}}}{2f^2(u)}} \,\rmd u,
\end{align}
where $a = e^{-\alpha} + e^{\alpha}$ and
\begin{align*}
    \bar{\Phi}_{s_1, y_1, s_2, y_2} := \bar{\Phi}\lrp{\frac{y_2 - y_1}{\sqrt{s_2 - s_1}}}, \quad  
    \phi_{s_1, y_1, s_2, y_2} := \phi\lrp{\frac{y_2 - y_1}{\sqrt{s_2 - s_1}}}, \quad y_1, y_2 \in \R, s_2 \geq s_1 \geq 0.
\end{align*}
Take now $y\downarrow b(s)$ in both \eqref{eq:pricing_formula} and \eqref{eq:pricing_formula_refined} to derive the free-boundary equation
\begin{align}\label{eq:free-boundary_eq}
    G(s, b(s)) &= c - \Es{\int_0^\infty \partial_t G\lrp{s + u, Y_u}\mathbbm{1}\lrp{Y_u \geq b(s + u)}\,\rmd u}{b(s)},
\end{align}
alongside its more explicit expression
\begin{align*}
    G&(s, b(s)) \\
    &= c - \int_s^\infty \frac{1}{f(u)}\lrp{c\bar{\Phi}_{s, b(s), u, b(u)} - \frac{(a + 2u)\lrp{(b(s) + cu)\bar{\Phi}_{s, b(s), u, b(u)} + \sqrt{u - s}\phi_{s, b(s), u, b(u)}}}{2f^2(u)}} \,\rmd u.
\end{align*}

It turns out that there exists a unique function $b$ that solves \eqref{eq:free-boundary_eq}, as we state in the next theorem. The proof of such an assertion follows from adapting the methodology used in \citet[Theorem 3.1]{Peskir_2005_American}, where the uniqueness of the solution of the free-boundary equation is addressed for an American put option with a geometric Brownian motion. % which had been an open problem for decades.

\begin{theorem}
    The integral equation \eqref{eq:free-boundary_eq} admits a unique solution among the class of continuous functions $\beta:\R_+\rightarrow\mathbb{R}$ of bounded variation. %and such that $\beta(s) > c$ for all $s\in\R_+$.
\end{theorem}

\begin{proof}
    % Intro, y -> W^beta(s,y) is twice continuously differentiable, and infinitesimal generator of Y acting on W^beta
	
    Suppose there exists a function $\beta:\R_+\rightarrow \mathbb{R}$ solving the integral equation \eqref{eq:free-boundary_eq}, and define $W^\beta$ as in \eqref{eq:pricing_formula}, but with $\beta$ instead of $b$. Since the integrand in \eqref{eq:pricing_formula_refined} is twice continuously differentiable with respect to $y$ and once with respect to $s$, we can use the Leibnitz rule to obtain $\partial_x W^\beta$, $\partial_{xx} W^\beta$, and $\partial_t W^\beta$ by differentiating within the integral symbol in \eqref{eq:pricing_formula_refined}, ensuring that these are continuous functions on $\R_+\times \R$. Besides, the following expression for $\bbL W^\beta$ can be easily computed from \eqref{eq:pricing_formula}:
    \begin{align*}
        \bbL W^\beta(s, y) = \partial_t G(t, y)\mathbbm{1}(y \geq \beta(s)).
    \end{align*}
	
    % Itô's formula for W^beta and G
	
    Define the sets 
    \begin{align*}
        \cC_\beta := \lrb{(s, y) \in \R_+\times\R : y < \beta(s)},\ \ 
        \cD_\beta := \lrb{(s, y) \in \R_+\times\R : y \geq \beta(s)}.
    \end{align*}
    Since $W^\beta \in C(\R_+\times\R)$, $\bbL W^\beta(s, y)$ is locally bounded, and $\beta$ is assumed to be continuous and of bounded variation, we can apply the \textit{iii-b} version of the Itô formula extension in Lemma A2 in \cite{DAuria_2020_discounted} (see \cite{Peskir_2005_change} for an original formulation) to obtain
    \begin{align}\label{eq:pricing_formula_W^beta}
        W^\beta(s, y) &= \Es{W^\beta(s + t, Y_t)}{y} - \Es{\int_0^t \partial_t G\lrp{s + u, Y_u}\mathbbm{1}\lrp{Y_u \geq \beta(s + u)}\,\rmd u}{y},
    \end{align}
    where the martingale term is canceled after taking $\Pr_y$-expectation and the local time term is missing due to the continuity of $\partial_x W^\beta$ on $\partial \cC_\beta$. In addition,
    \begin{align}\label{eq:G_Ito}
        G(s, y) &= \Es{G(s + t, Y_t)}{y} - \Es{\int_0^t \partial_t G\lrp{s + u, Y_u}\,\rmd u}{y}.
    \end{align}
	
    % W^beta = G on D_beta!!
	
    Due to the law of the iterated logarithm, and recalling \eqref{eq:gain} and \eqref{eq:G_t}, we obtain
    \begin{align}\label{eq:G-lim}
        \lim_{u\rightarrow\infty}G(s + u, Y_u) &= c 
    \end{align}
    and 
    \begin{align}\label{eq:G_t-lim}
        \lim_{u\rightarrow\infty}\partial_t G\lrp{s + u, Y_u} &= 0
    \end{align}
    $\Pr_y$-a.s. for all $y\in\R$. We get the following from the fact that $W^\beta$ satisfies \eqref{eq:pricing_formula} with $\beta$ instead of $b$, along with the dominated convergence theorem and \eqref{eq:G_t-lim}:
    \begin{align}\label{eq:W-lim}
        \lim_{u\rightarrow\infty}W^\beta(s + u, Y_u) = c
    \end{align} 
    $\Pr_y$-a.s. for all $y\in\R$. Consider the first hitting time $\sigma_{\cC_\beta}$ into $\cC_\beta$, fix $(s, y)\in\cD_\beta$. From \eqref{eq:G-lim} and \eqref{eq:W-lim}, and the fact that $W^\beta(s, \beta(s)) = G(s, \beta(s))$ for all $s\in\R_+$ (as $\beta$ solves \eqref{eq:free-boundary_eq}), it follows that $W^\beta\big(s + \sigma_{\cC_\beta}, Y_{\sigma_{\cC_\beta}}\big) = G\big(s + \sigma_{\cC_\beta}, Y_{\sigma_{\cC_\beta}}\big)\ \Pr_y$-a.s. for all $y\in\R$ and all $s\in\R_+$. Relying on this last identity along with the fact that $\Pr_y(Y_u \geq \beta(t + s)) = 1$ for all $0 \leq u \leq \sigma_{\cC_\beta}$, and using \eqref{eq:pricing_formula_W^beta} and \eqref{eq:G_Ito}, we obtain
    \begin{align*}
        W^\beta(s, y) &= \Es{W^\beta\big(s + \sigma_{\cC_\beta}, Y_{\sigma_{\cC_\beta}}\big)}{y} - \Es{\int_0^{\sigma_{\cC_\beta}} \partial_t G\lrp{s + u, Y_u}\,\rmd u}{y} \\
        &= \Es{G^\beta\big(s + \sigma_{\cC_\beta}, Y_{\sigma_{\cC_\beta}}\big)}{y} - \Es{\int_0^{\sigma_{\cC_\beta}} \partial_t G\lrp{s + u, Y_u}\,\rmd u}{y} \\
	& = G(s, y),
    \end{align*}
    which proves that $W^\beta = G$ on $\cD_\beta$.
	
    % W^beta <= W!!
	
    Define now the first hitting time $\sigma_{\cD_\beta}$ into $\cD_\beta$. Note that either $\sigma_{\cD_\beta} = 0$ for $(s, y)\in\cD_\beta$, on which $W^\beta =  G$, or $Y_u < \beta(s + u)$ for all $0\leq u < \sigma_{\cD_\beta}$. We derive from  \eqref{eq:pricing_formula_W^beta} that
    \begin{align*}
        W^\beta(s, y) &= \Es{W^\beta\lrp{s + \sigma_{\cD_\beta}, Y_{\sigma_{\cD_\beta}}}}{y} = \Es{G\lrp{s + \sigma_{\cD_\beta}, Y_{\sigma_{\cD_\beta}}}}{y}
    \end{align*}
    for all $(s, y)\in\R_+\times\R$, which, after recalling the definition of $W$ in \eqref{eq:OSP_BM}, proves that $W^\beta \leq W$.
	
    % b >= beta!!
	
    Take $(s, y)\in \cD_\beta\cap\cD$ and consider the first hitting time $\sigma_\cC$ into the continuation set $\cC$. Since \mbox{$W = G$} on $\cD$ and $W^\beta = G$ on $\cD_\beta$, by relying on \eqref{eq:pricing_formula}, \eqref{eq:pricing_formula_W^beta}, and the fact that \mbox{$\Pr_y\lrp{Y_u \geq b(s + u)} = 1$} for all $0 \leq u < \sigma_\cC$, we get
    \begin{align*}
        \Es{W\lrp{s + \sigma_\cC, Y_{\sigma_\cC}}}{y} &= G(s, y) + \Es{\int_0^{\sigma_\cC} \partial_t G\lrp{s + u, Y_u}\,\rmd u}{y}, \\
        \Es{W^\beta\lrp{s + \sigma_\cC, Y_{\sigma_\cC}}}{y} &= G(s, y) + \Es{\int_0^{\sigma_\cC} \partial_t G\lrp{s + u, Y_u}\mathbbm{1}\lrp{Y_u \geq \beta(s + u)}\,\rmd u}{y}.
    \end{align*}
    After recalling that $W^\beta \leq W$, we can merge the two previous equalities into
    \begin{align*}
        \Es{\int_0^{\sigma_\cC} \partial_t G\lrp{s + u, Y_u}\mathbbm{1}\lrp{Y_u \geq \beta(s + u)}\,\rmd u}{y} &\leq \Es{\int_0^{\sigma_\cC} \partial_t G\lrp{s + u, Y_u}\,\rmd u}{y},
    \end{align*}
    which, alongside the fact that $\partial_t G(s, y) < 0$ for all $(s, y)\in\cD$ (otherwise we get from \eqref{eq:pricing_formula_aux} that the first exit time from a ball around $(s, y)$ small enough will yield a better strategy than stopping immediately) and the continuity of $\beta$, implies that $b \geq \beta$. 
	
    % b = beta!!!!!  At last!
	
    Suppose that there exists a point $s\in\R_+$ such that $b(s) > \beta(s)$ and fix $y\in(\beta(s), b(s))$. Consider the stopping time $\sigma^* = \sigma^*(s, y)$ and plug it into both \eqref{eq:pricing_formula} and \eqref{eq:pricing_formula_W^beta} to obtain
    \begin{align*}
        \Es{W^\beta\lrp{s + \sigma^*, Y_{\sigma^*}}}{y} &= \Es{G\lrp{s + \sigma^*, Y_{\sigma^*}}}{y}\\
        &= W^\beta(s, y) + \Es{\int_0^{\sigma^*} \partial_t G\lrp{s + u, Y_u}\mathbbm{1}\lrp{Y_u \geq \beta(s + u)}\,\rmd u}{y}
    \end{align*}
    and
    \begin{align*}
        \Es{W\lrp{s + \sigma^*, Y_{\sigma^*}}}{y} = \Es{G\lrp{s + \sigma^*, Y_{\sigma^*}}}{y} = W(s, y).
    \end{align*} 
    Thus, since $W^\beta \leq W$, we get
    \begin{align*}
        \Es{\int_0^{\sigma^*} \partial_t G\lrp{s + u, Y_u}\mathbbm{1}\lrp{Y_u \geq \beta(s + u)}\,\rmd u}{y} \geq 0.
    \end{align*}
    Using the fact that $y < b(s)$, the continuity of $b$, and the time-continuity of the process $Y$, we can state that $\sigma^* > 0\ \Pr_y$-a.s. Therefore, since $\partial_t G(s, y) < 0$ for all $(s, y)\in\cD_\beta$ (the same arguments used to prove that $\partial_t G < 0$ in $\cD$ lead to this conclusion) the previous inequality can only stand if $\mathbbm{1}\lrp{Y_u \geq \beta(s + u)} = 0$ for all $0\leq u\leq \sigma^*$, meaning that $b(s + u) \leq \beta(s + u)$ in the same interval, which contradicts the assumption $b(s) > \beta(s)$ due to the continuity of both $b$ and $\beta$. 
\end{proof}

%-------------------------------------------------%
\section{Solution of the original problem and some extensions}\label{sec:sol_original}
%-------------------------------------------------%

Recall that the OSPs \eqref{eq:OSP_BM} and \eqref{eq:OSP_OUB} are equivalent, meaning that the value functions and the OSTs of both problems are linked through a homeomorphic transformation. Details on how to actually translate one problem into the other were given in Proposition \ref{pr:OSP_equiv}. It then follows that the stopping time $\tau^*(t, x)$ defined in \eqref{eq:OST_transform} is optimal for \eqref{eq:OSP_BM} and it admits the following alternative representation under~$\Pr_x$:
\begin{align}\label{eq:OST_OSB}
    \tau^*(t, x) = \inf\lrb{u\geq 0 : X_{t + u} \geq \beta(t + u)},\quad \beta(t) = \frac{z}{c_z}G_{c_z}\lrp{s, b(s)},
\end{align}
where $\beta$ is the OSB associated to \eqref{eq:OSP_OUB}, and $s = \upsilon(t)$ and $c_z = \eta(z)$ are defined in Proposition \ref{pr:OSP_equiv}. We can obtain both $V$ and $\beta$ without requiring the computation of $W$ and $b$. Indeed, consider the infinitesimal generator of $\lrb{\lrp{t, X_t}}_{t \in [0, 1]}$, $\bbL_X$, and set $y = \eta(x)$, $s_\varepsilon = s + \varepsilon$, and $t_\varepsilon = \upsilon^{-1}(s_\varepsilon)$ for $\varepsilon\in\R$. By means of \eqref{eq:value_equiv} and the chain rule, we get that
\begin{align*}
    \frac{z}{c_z}\lrp{\bbL W_{c_z}}(s, y) :=&\; \lim_{\varepsilon\rightarrow 0}\varepsilon^{-1}\lrp{\Es{\frac{z}{c_z}W_{c_z}\lrp{s_\varepsilon, Y_\varepsilon}}{y} - \frac{z}{c_z}W_{c_z}(s, y)} \\
    =&\; \lim_{\varepsilon\rightarrow 0}\varepsilon^{-1}\lrp{\Es{V(t_\varepsilon, X_{t_\varepsilon})}{t, x} - V(t, x)} \\
    =&\; \lrp{\bbL_X V}(t, x)\big(\upsilon^{-1}\big)'(s).
\end{align*}
Hence, after multiplying both sides of \eqref{eq:pricing_formula_aux} by $z/c_z$, integrating with respect to $\upsilon^{-1}(u)$ instead of $u$, and recalling that $\bbL_X V(t, x) = 0$ for all $x\leq \beta(t)$ and $V(t, x) = x$ for all $x\geq\beta(t)$, we get the pricing formula
\begin{align}
    V(t, x) &= z - \Es{\int_0^{1-t} (\bbL_X V)(t + u, X_{t + u})\,\rmd u}{t, x} \nonumber \\
    &= z - \Es{\int_0^{1-t} \mu(t + u, X_{t + u})\mathbbm{1}(X_{t + u} \geq \beta(t + u))\,\rmd u}{t, x}. \label{eq:pricing_formula_original}
\end{align}
In the same fashion we obtained \eqref{eq:pricing_formula_refined}, we can take advantage of the linearity of $x\mapsto\mu(t, x)$ and the Gaussian marginal distributions of $X$ to come up with the following refined version of \eqref{eq:pricing_formula_original}:
\begin{align}\label{eq:pricing_formula_original_refined}
    V(t, x) = z - \int_t^1K(t, x, u, \beta(u))\,\rmd u,
\end{align}
where, for $x_1, x_2 \in\R$ and $0\leq t_1 \leq t_2 \leq 1$,
\begin{align}\label{eq:kernel}
    K(t_1, x_1, t_2, x_2) := \alpha\frac{z\wt{\Phi}_{t_1, x_1, t_2, x_2} - \cosh(\alpha(1 - t_2))(m_{t_2}(t_1, x_1)\wt{\Phi}_{t_1, x_1, t_2, x_2} + v_{t_2}(t_1)\wt{\phi}_{t_1, x_1, t_2, x_2})}{\sinh(\alpha(1 - t_2))},
\end{align}
with
\begin{align*}
    \wt{\Phi}_{t_1, x_1, t_2, x_2} := \bar{\Phi}\lrp{\frac{x_2 - m_{t_2}(t_1, x_1)}{v_{t_2}(t_1)}}, \quad  
    \wt{\phi}_{t_1, x_1, t_2, x_2} := \phi\lrp{\frac{x_2 - m_{t_2}(t_1, x_1)}{v_{t_2}(t_1)}}
\end{align*}
and 
\begin{align*}
    m_{t_2}(t_1, x_1) &:= \Es{X_{t_2}}{t_1, x_1} = \frac{x_1\sinh(\alpha(1 - t_2)) + z\sinh(\alpha (t_2 - t_1))}{\sinh(\alpha(1 - t_1))}, \\
    v_{t_2}(t_1) &:= \sqrt{\Vs{X_{t_2}}{t_1}} = \sqrt{\frac{\gamma^2}{\alpha}\frac{\sinh(\alpha(1 - t_2))\sinh(\alpha (t_2 - t_1))}{\sinh(\alpha(1 - t_1))}}.
\end{align*}
Consequently, by taking $x\downarrow\beta(t)$ in \eqref{eq:pricing_formula_original} (or by directly transforming \eqref{eq:free-boundary_eq} in the same way we obtained \eqref{eq:pricing_formula_original} from \eqref{eq:pricing_formula}), we get the free-boundary equation
\begin{align*}
    \beta(t) &= z - \Es{\int_0^{1-t} (\bbL_X V)(t + u, X_{t + u})\,\rmd u}{t, \beta(t)} \\
    &= z - \Es{\int_0^{1-t} \mu(t + u, X_{t + u})\mathbbm{1}(X_{t + u} \geq \beta(t + u))\,\rmd u}{t, \beta(t)},
\end{align*}
which may also be expressed as
\begin{align}\label{eq:free-boundary_eq_original_refined}
    \beta(t) = z - \int_t^1 K(t, \beta(t), u, \beta(u))\,\rmd u.
\end{align}

The next three remarks broaden the scope of applicability of the OUB as the underlying model in \eqref{eq:OSP_OUB}. In particular, the two first reveal that setting the terminal time to $1$ and the pulling level (coming from the asymptotic mean of the OU process underneath) to $0$ does not take a toll on generality, while the last one shows that the OSP for the BB arises as a limit case when $\alpha\to0$.

\begin{remark}[OUB with a general pulling level]\label{rmk:general_pulling_level}
    Let $\smash{\wt{X}^\theta = \big\{\wt{X}_t^\theta\big\}_{t \in [0, 1]}}$ be an OU process satisfying the SDE $\rmd \wt{X}_t^\theta = \alpha(\wt{X}_t^\theta - \theta)\,\rmd t + \gamma\,\rmd B_t$. That is, $X^{\theta, z}$ is pulled towards $\theta$ with a time-dependent strength dictated by $\alpha$. Denote by $\smash{X^{\theta, z} = \big\{X_t^{\theta, z}\big\}_{t \in [0, 1]}}$ the OUB process built on top of $\wt{X}^\theta$ and such that $X_1^{\theta,z} = z$. It is easy to check that 
    %(PROVE OF "IT IS EASy TO CHECK")
    %$\wt{X}_t^\theta = \wt{X}_t + \theta(1 - e^{-\alpha t})$, from which follows the relation $X_t^\theta = X_t + \theta\left(1 - \frac{\sinh(\alpha (1 - t)) + \sinh(\alpha t)}{\sinh(\alpha)}\right)$. Then, since $X^\theta$ admits (see\cite{Barczy2013Sample}) the integral representation
    %\begin{align*}
	%X_t = x\frac{\sinh(\alpha (1 - t))}{\sinh(\alpha)} + z\frac{\sinh(\alpha t)}{\sinh(\alpha)}  + \gamma\int_{0}^t \frac{\sinh(\alpha (1 - t))}{\sinh(\alpha (1 - u))}\,\mathrm{d}B_u,
	%\end{align*}
	$X^{\theta, z} = X^{0, z - \theta} + \theta$, whenever $X_0^{0, z - \theta} = X_0^{\theta, z} - \theta$. Denote by $V^{\theta, z}$ and $\beta^{\theta, z}$ the value function and the OSB related to the OSP \eqref{eq:OSP_OUB} with $X$ replaced by $X^{\theta, z}$. Then $V^{\theta, z}(t, x) = V^{0, z - \theta}(t, x - \theta) + \theta$ and $b^{\theta, z}(t) =  b^{0, z - \theta}(t) + \theta$.
\end{remark}

\begin{remark}[OUB with a general horizon]
    Denote by $\smash{X^{\alpha, \gamma, T} = \big\{X_t^{\alpha, \gamma, T}\big\}_{t \in [0, T]}}$ an OUB with slope $\alpha$, volatility $\gamma$, and horizon $T$. Likewise, let $V^{\alpha, \gamma, T}$ and $\beta^{\alpha, \gamma, T}$ be the corresponding value function and the OSB. By relying on the scaling property of a Brownian motion, one can easily verify that $X_t^{\alpha r, \gamma, T} = X_{rt}^{\alpha, \gamma r^{-1/2}, rT}$ $\Pr_x$-a.s. for any $r > 0 $. Consequently, $V^{\alpha r, \gamma, T}(t, x) = V^{\alpha, \gamma r^{-1/2}, rT}(rt, x)$ and $\beta^{\alpha r, \gamma, T}(t) = \beta^{\alpha, \gamma r^{-1/2}, rT}(rt)$. Thereby, by taking $r = 1/T$, one can derive $V^{\alpha, \gamma, T}$ and $\beta^{\alpha, \gamma, T}$ for any set of values $\alpha$, $\gamma$, and $T$ from the solution of the OSP in \eqref{eq:OSP_OUB}.
\end{remark}

\begin{remark}[BB from an OUB]
    To emphasize the dependence on $\alpha$, denote by $X(\alpha)$, $V_\alpha$, and $\beta_\alpha$, respectively, the OUB solving \eqref{eq:OUB_SDE}, the value function in \eqref{eq:value}, and the corresponding OSB. The process $X_t(\alpha)$ has the following integral representation under $\Pr_x$ \citep{Barczy_2013_sample}:
    \begin{align*}
	X_{t} = x\frac{\sinh(\alpha (1 - t))}{\sinh(\alpha)} + z\frac{\sinh(\alpha t)}{\sinh(\alpha)} + \sigma\int_{0}^t \frac{\sinh(\alpha (1 - t))}{\sinh(\alpha (1 - u))}\,\rmd B_u,
	\end{align*}
	from where we can conclude, after taking $\alpha\rightarrow 0$ and using the dominated convergence theorem, that $X_t(\alpha) \rightarrow \wt{X}_t$ $\Pr_x$-a.s. for all $t\in[0, 1)$, where $\wt{X}$ is a BB process with final value $\tilde{X}_1 = z$. Then, by applying Theorem 5 from \cite{Coquet_2007_convergence} we have that $V_\alpha \rightarrow \wt{V}$, and hence $\beta_\alpha \rightarrow \wt{\beta}$, as $\alpha\rightarrow 0$, where $\wt{V}$ and $\wt{\beta}$ are the value function and the OSB related to $\wt{X}$.
\end{remark}

\begin{remark}[Time-dependent gain function]
    Note that the same methodology used to obtain the solution of \eqref{eq:OSP_OUB}, that is, the time-space transformation of the GMB into a BM, might be extended to cover a wider class of time-dependent gain functions. Indeed, let
    \begin{align*}
        \wt{V}(t, x) := \sup_{\tau \leq 1 - t}\Es{\wt{F}\lrp{t + \tau, X_{t + \tau}}}{t, x},
    \end{align*}
    for $\wt{F}:[0, T]\times \R \mapsto \R$. Then, using the same notation and arguing as in Proposition \ref{pr:OSP_equiv}, one can get that
    \begin{align*}
        \wt{V}(t, x) := \wt{W}_{c_{z}}(s, y),
    \end{align*}
    for
    \begin{align*}
        \wt{W}_c(s, y) := \sup_{\sigma}\Es{\wt{G}_c(s + \sigma, Y_\sigma)}{y},\quad \wt{G}_c(s, y) := \wt{F}\lrp{\nu^{-1}(s), \frac{z}{c_z}G_{c_z}(s, y)}.
    \end{align*}
    Therefore, solving $\wt{V}$ and $\wt{W}$ are equivalent problems as happens with $V$ and $W$. As long as $\wt{F}$ is linear in its space coordinate, the applicability of the methods used in Section \ref{sec:sol_reformulated} to solve $\wt{W}$ might come straightforwardly provided smoothness of its time coordinate, namely, among others, partial derivability, boundedness, and limit as time diverges. The time-smoothness of $\wt{F}$ is inherited from that of $\wt{G}$ and the change of time $\upsilon$.
\end{remark}

\begin{remark}[Gauss--Markov bridges]
    The methodology we used to obtain the solution to the OSP \eqref{eq:OSP_OUB} is essentially based on the time-space equivalence \eqref{eq:OUB_to_BM}, which allows us to work in the simpler Brownian motion scenario. Such types of representations are not exclusive to the OUB, but they are shared by processes that result from conditioning Gauss--Markov processes to hit a deterministic terminal point (see, e.g., \cite{Mehr_1965_certain}, \cite{Borisov_1983_criterion}, and \cite{Buonocore_2013_some}). Therefore, the OSPs of processes from this wider class of bridges can be addressed by following a Brownian-motion representation approach similar to the one used in this paper to solve the OUB case. The recent work of \cite{Azze_2022_GMB} further develops this idea.
\end{remark}

%-------------------------------------------------%
\section{Numerical results}\label{sec:numerical_results}
%-------------------------------------------------%

The free-boundary equation \eqref{eq:free-boundary_eq_original_refined} does not admit a closed-form solution and thus numerical procedures are needed to compute an approximate boundary. 
By exploiting the fact that the OSB at a given time $t$ depends only on its shape from $t$ up to the horizon, one can discretize the integral in \eqref{eq:free-boundary_eq_original_refined} by means of a right Riemann sum and, since the terminal value $\beta(1)$ is known, the entire boundary can be computed in a backward form. This method of backward induction is detailed in \citet[Chapter 8]{Detemple_2005_American-style} and examples of its implementation can be found, e.g., in \cite{Pedersen_2002_onnonlinear}.
Another approach to solve \eqref{eq:free-boundary_eq_original_refined} is by using Picard iterations, that is, by treating \eqref{eq:free-boundary_eq_original_refined} as a fixed-point problem in which the entire boundary is updated in each step. The works of \cite{Detemple_2020_value} and \cite{DeAngelis_2020_optimal} use this approach to solve the associated Volterra-type integral equation characterizing the OSB. To the best of our knowledge, when it comes to non-linear integral equations arisen from OSPs, the convergence of both the Picard scheme and the backward induction technique are numerically checked rather than formally proved. Therefore, we chose to use the Picard scheme since empirical tests suggested a faster convergence rate while keeping a similar accuracy compared to the backward induction approach.

Define a partition of $[0, 1]$, namely, $0 = t_0 < t_1 < \cdots < t_N = 1$ for $N\in\mathbb{N}$. Given that $\beta(1) = z$, we will initialize the Picard iterations by starting with the constant boundary $\beta^{(0)}:[0, 1] \to \R$ with $\beta^{(0)}\equiv z$. The updating mechanism that generates subsequent boundaries is laid down in the following formula, which comes after discretizing the integral in \eqref{eq:free-boundary_eq_original_refined} by using a right Riemann sum:
\begin{align*}
    \beta_i^{(k)} = z - \sum_{j = i}^{N - 2} K\lrp{t_i, \beta_i^{(k - 1)}, t_{j + 1}, \beta_{j + 1}^{(k - 1)}}(t_{j + 1} - t_j), \quad k = 1, 2, \dots
\end{align*}
We neglect the $(N-1)$-addend and allow the sum to run only until $N - 2$ since $K(t, x, 1, z)$ is not well defined, and therefore the last integral piece cannot be included in the right Riemann sum. As the overall integral is finite, the last piece vanishes as $t_{N - 1}$ approaches $1$. 

We chose to stop the fixed-point Picard algorithm after the $m$-th iteration if $m = \min\big\{k > 0: \max_{i = 1, \dots, N}| \beta_i^{k-1} - \beta_i^k| < \varepsilon\big\}$ for $\varepsilon=10^{-4}$. Empirical evidence suggested that the best performance of the algorithm was achieved when using a non-uniform mesh whose distances $t_i - t_{i-1}$ smoothly decrease as $i$ increases. In our computations, we used the logarithmically-spaced partition $t_i = \ln\lrp{1 + i(e - 1)/N}$, where $N = 500$ unless is otherwise specified.

Figures \ref{fig:alpha_change}, \ref{fig:gamma_change}, and \ref{fig:pinning_change} reveal how the OSB's shape is affected by different sets of values for the slope $\alpha$, the volatility $\gamma$, and the anchor point $z$.

\begin{figure}[!ht]
	\centering
	\begin{subfigure}[b]{0.32\textwidth}
		\includegraphics[width = \textwidth]{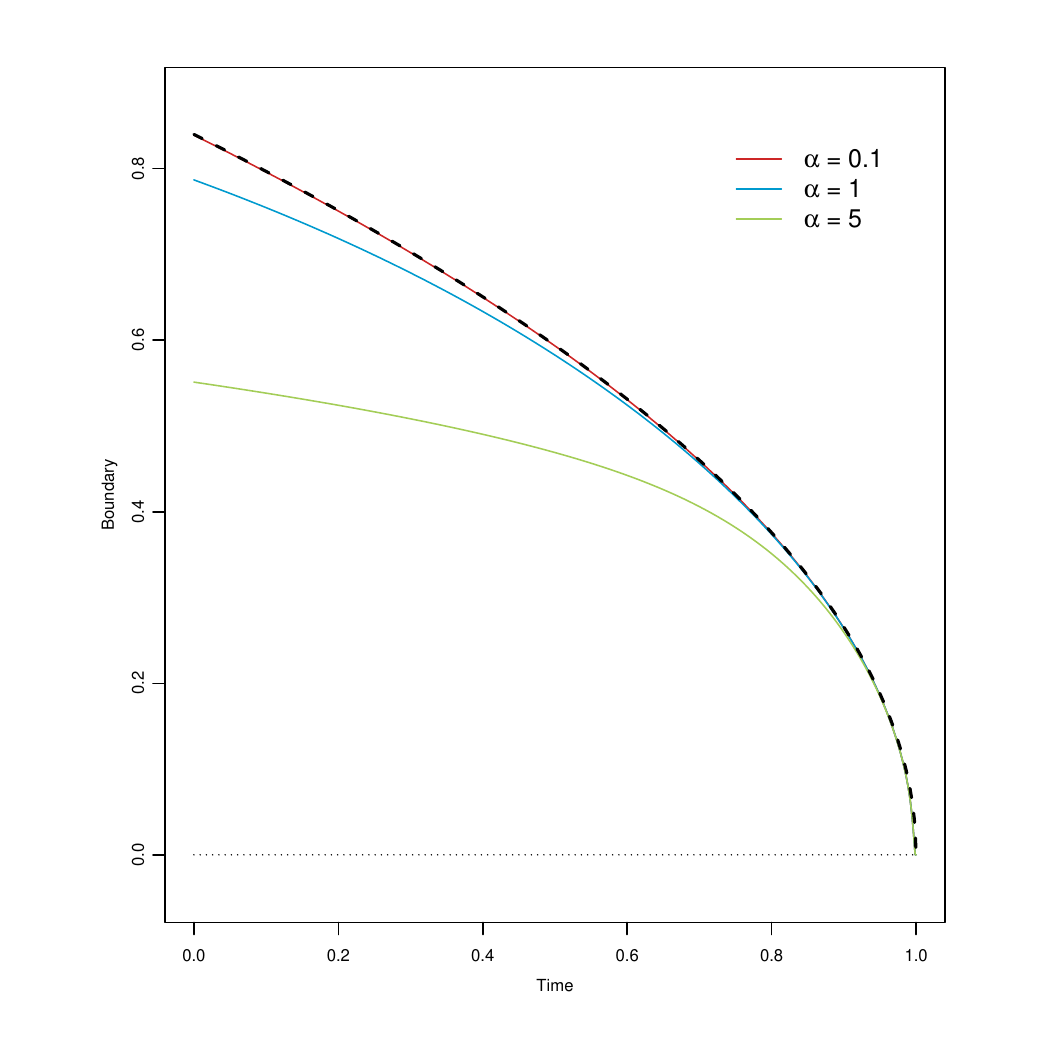}
		\subcaption{$z = 0, \gamma = 1$}
	\end{subfigure}
	\begin{subfigure}[b]{0.32\textwidth}
		\includegraphics[width = \textwidth]{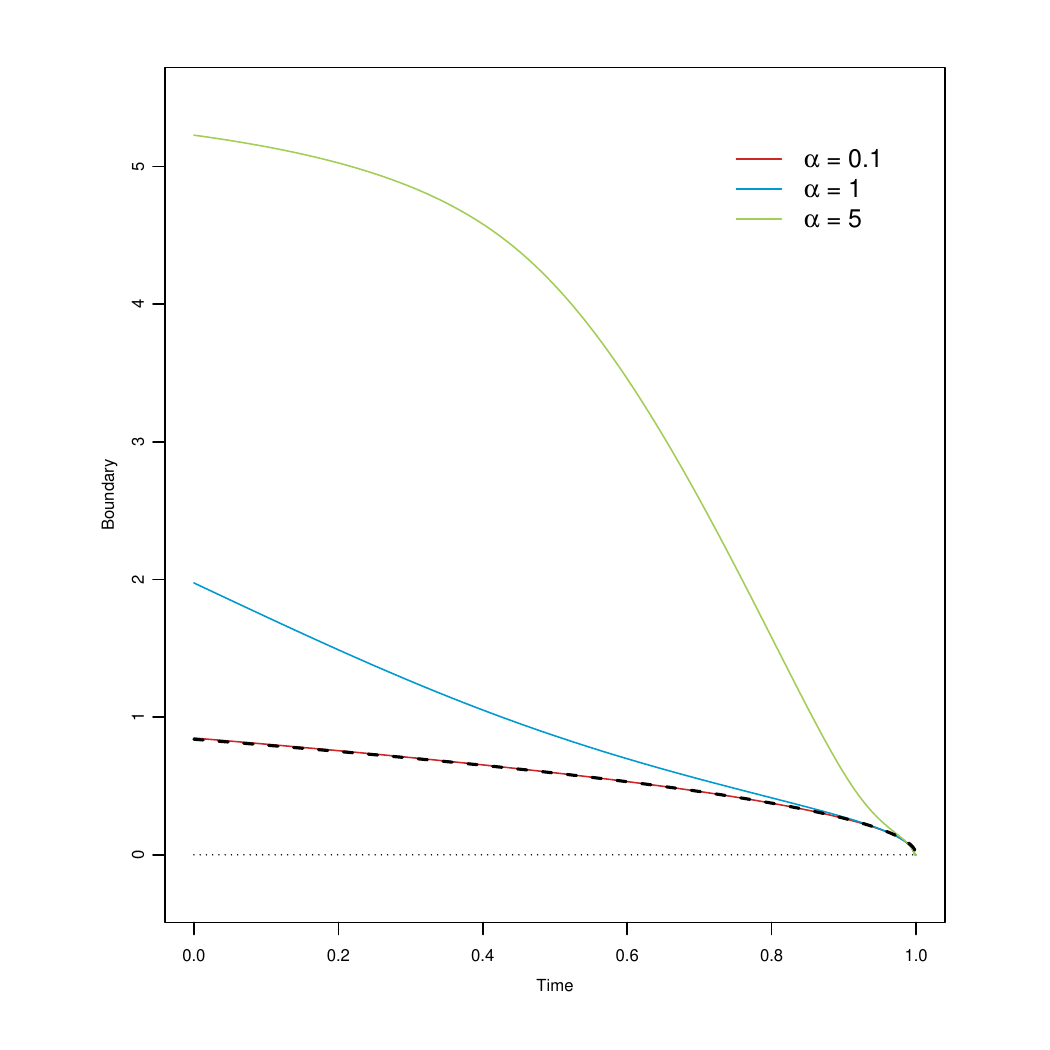}
		\subcaption{$z = -5, \gamma = 1$}
	\end{subfigure}
	\begin{subfigure}[b]{0.32\textwidth}
		\includegraphics[width = \textwidth]{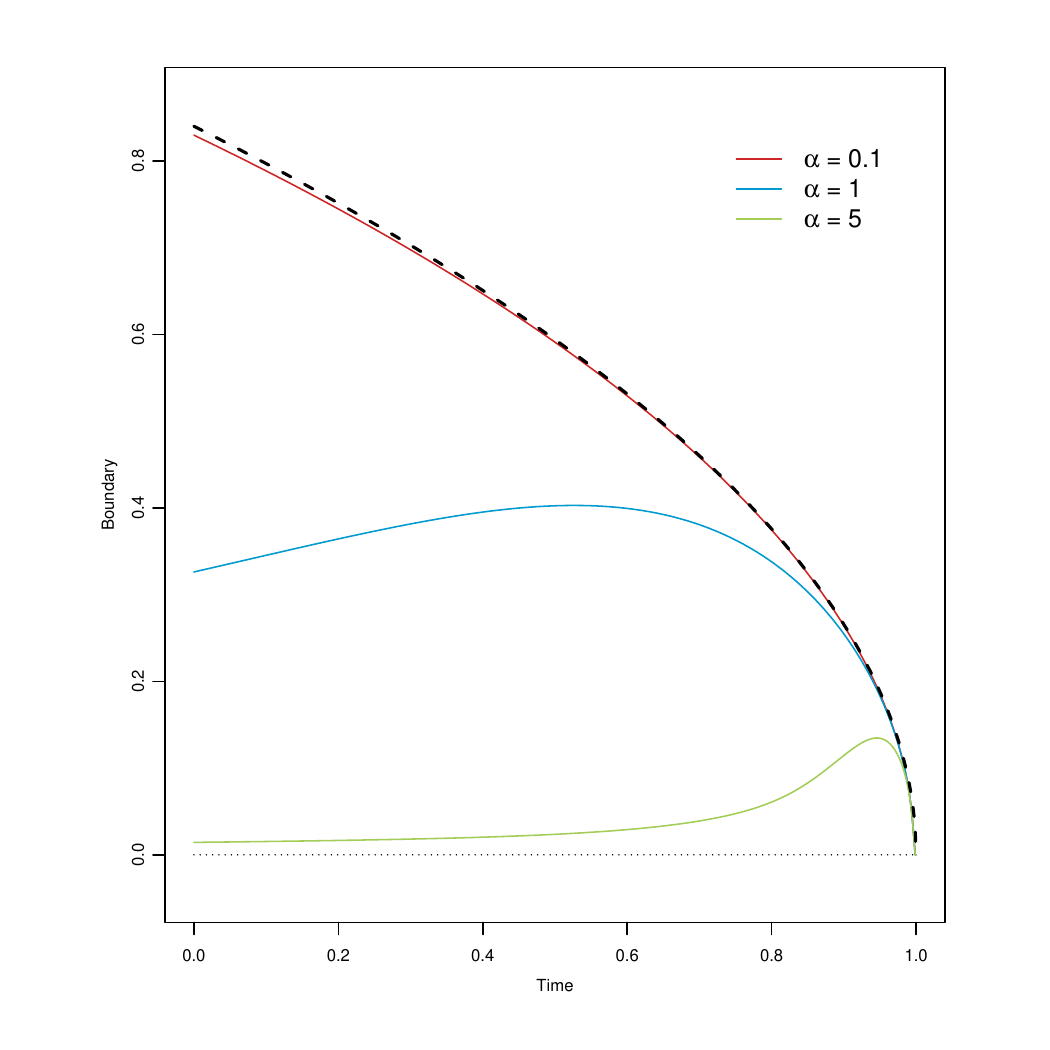}
		\subcaption{$z = 5, \gamma = 1$}
	\end{subfigure}
	\vspace{6PT}

	\caption{\small Optimal stopping boundary estimation for different values of $\alpha$. The boundary is pulled towards $0$ with a strength that increases as both $|\alpha|$ (values of $\alpha$ with equal absolute values yield the same boundary) and the residual time to the horizon $1 - t$ increases. As $\alpha \rightarrow 0$, the boundary estimation is shown to converge towards the OSB of a BB (dashed line), which is known to be $z + L\sqrt{1 - t}$, for $L\approx 0.8399$.}
	\label{fig:alpha_change}
\end{figure}

\begin{figure}[!ht]
	\centering
	\begin{subfigure}[b]{0.32\textwidth}
		\includegraphics[width = \textwidth]{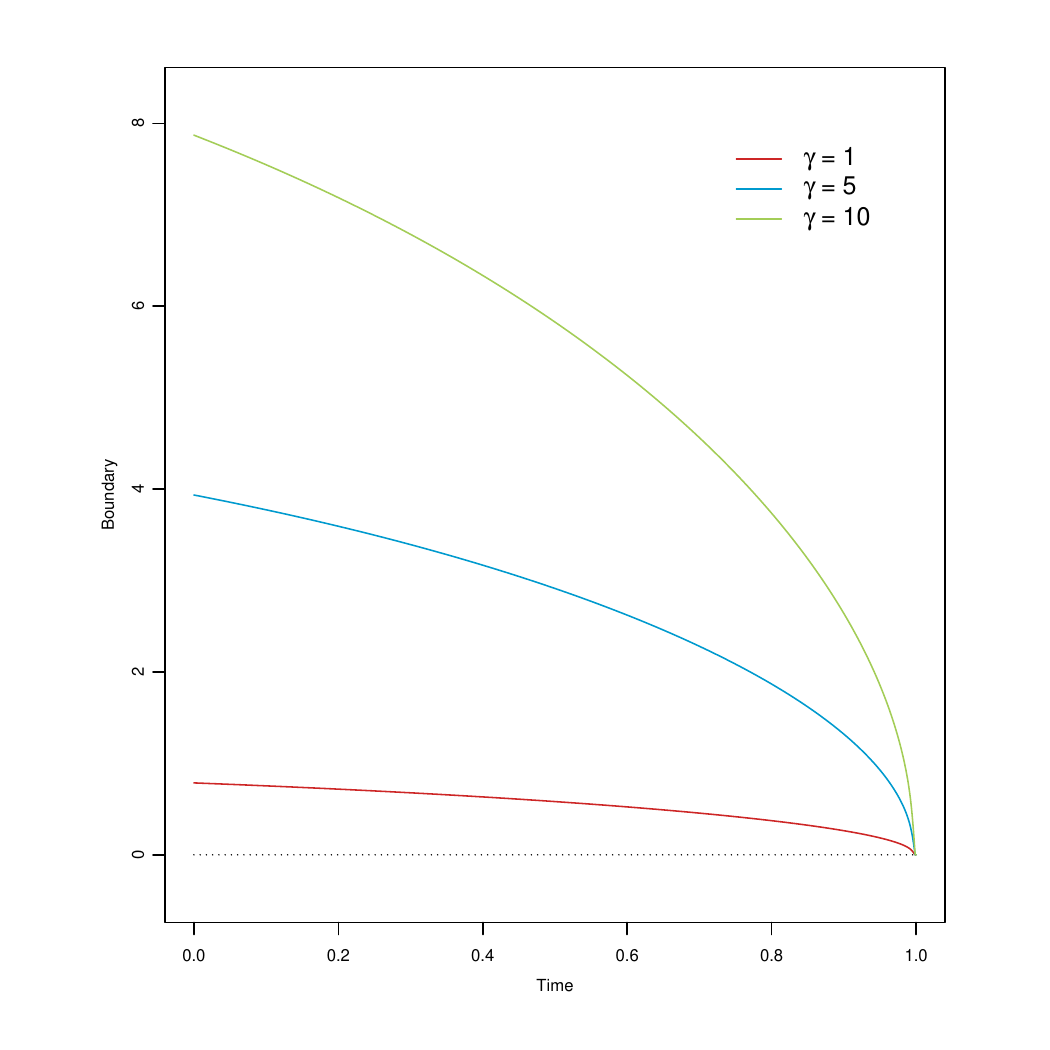}
		\subcaption{$z = 0, \alpha = 1$}
	\end{subfigure}
	\begin{subfigure}[b]{0.32\textwidth}
		\includegraphics[width = \textwidth]{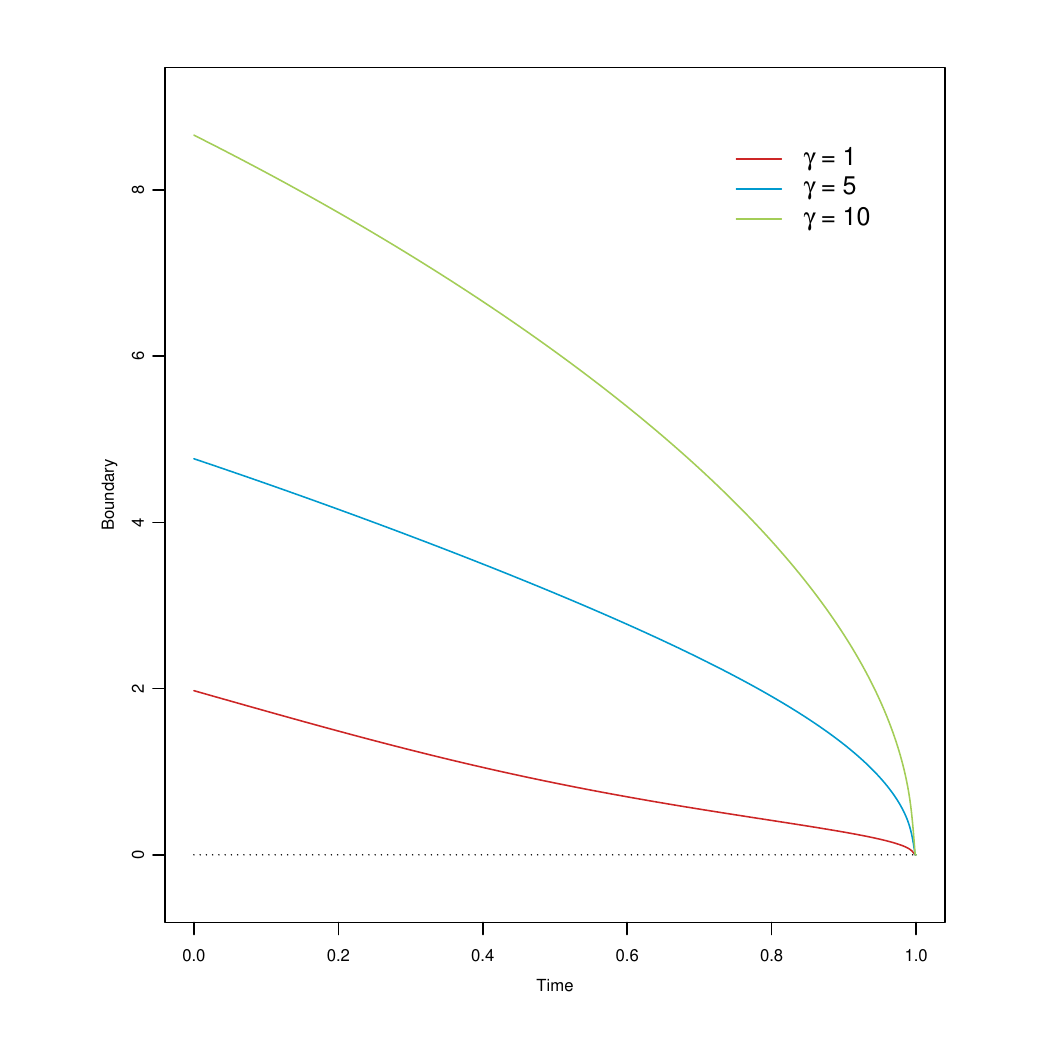}
		\subcaption{$z = -5, \alpha = 1$}
	\end{subfigure}
	\begin{subfigure}[b]{0.32\textwidth}
		\includegraphics[width = \textwidth]{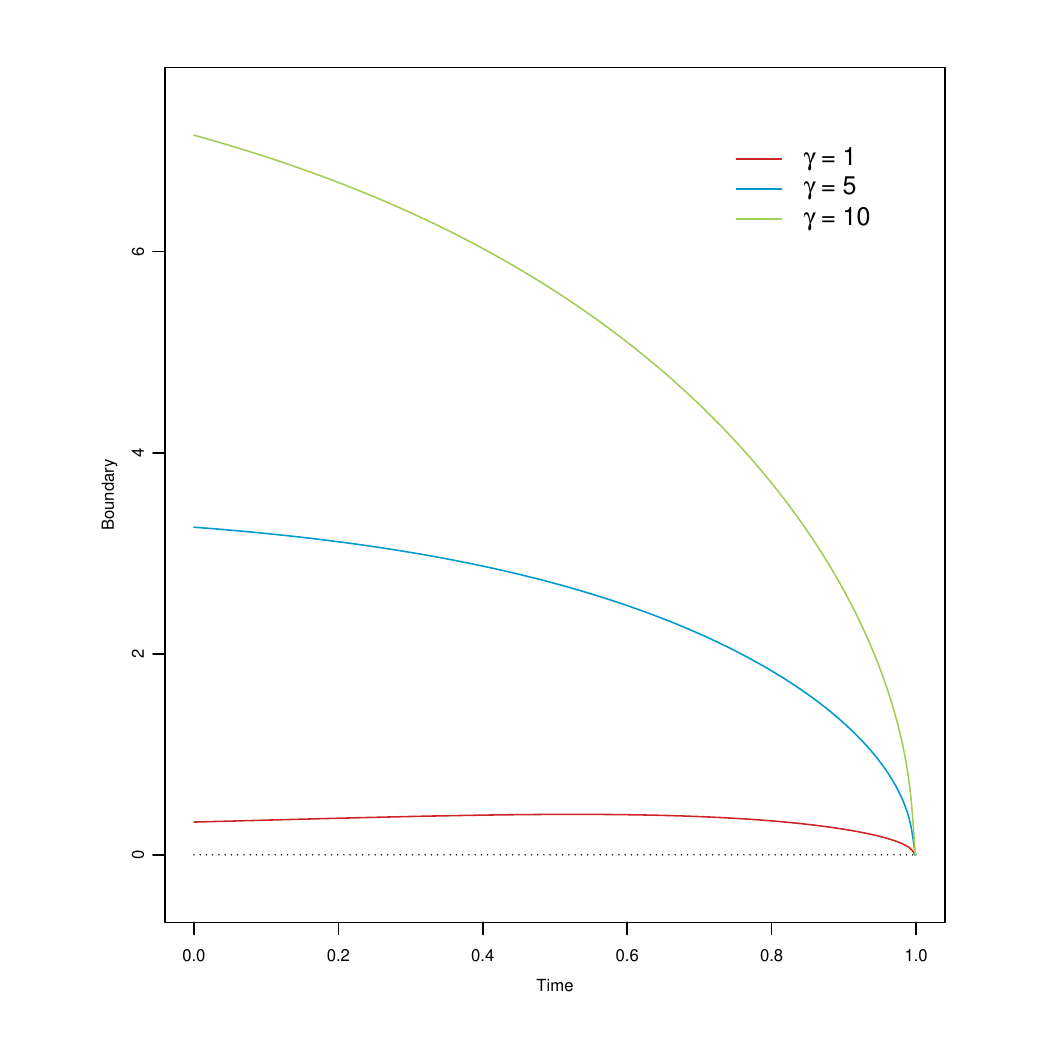}
		\subcaption{$z = 5, \alpha = 1$}
	\end{subfigure}
	\vspace{6PT}

	\caption{\small Optimal stopping boundary estimation for different values of $\gamma$. The boundary exhibits an increasing proportional relationship with respect to $\gamma$.}
	\label{fig:gamma_change}
\end{figure}

\begin{figure}[!ht]
	\centering
	\begin{subfigure}[b]{0.32\textwidth}
		\includegraphics[width = \textwidth]{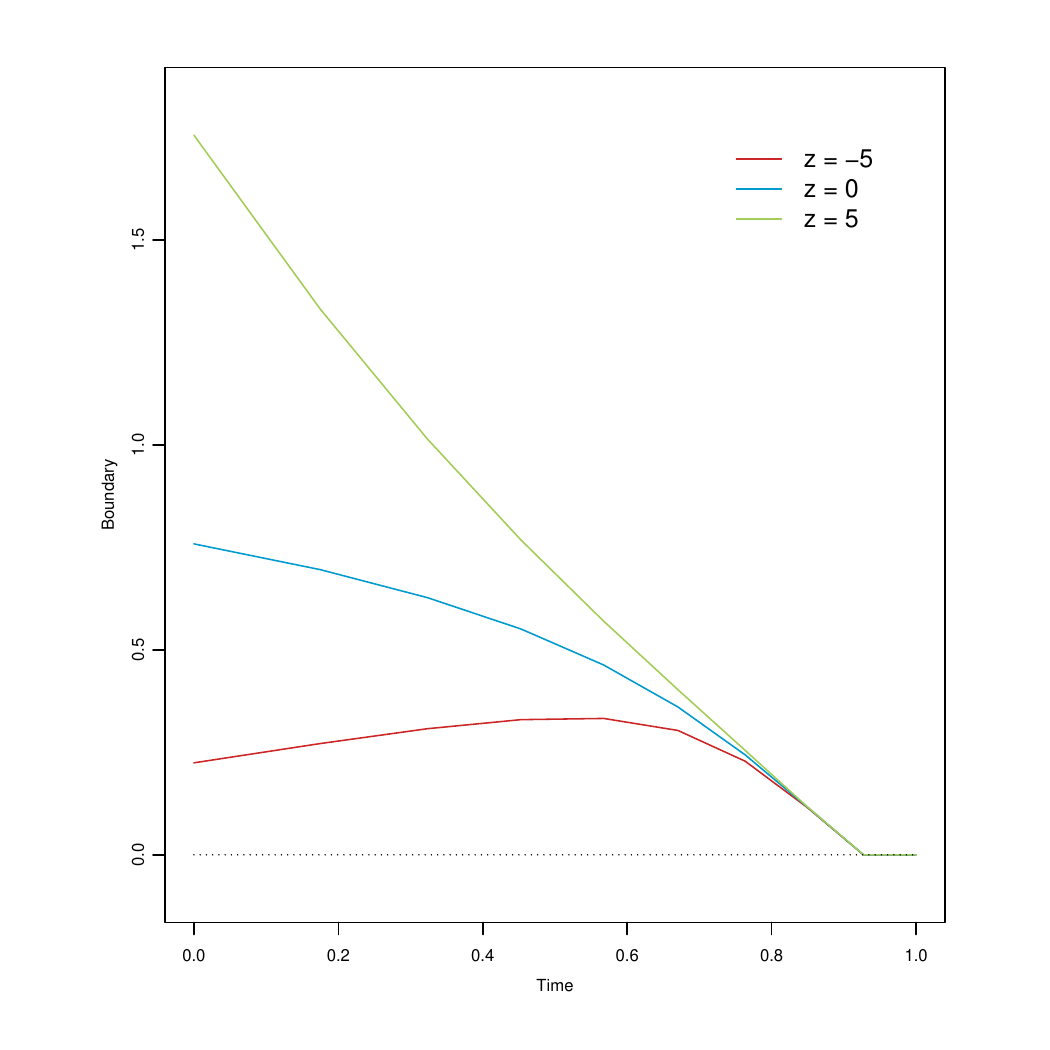}
		\subcaption{$\alpha = \gamma = 1, N = 10$}
	\end{subfigure}
	\begin{subfigure}[b]{0.32\textwidth}
		\includegraphics[width = \textwidth]{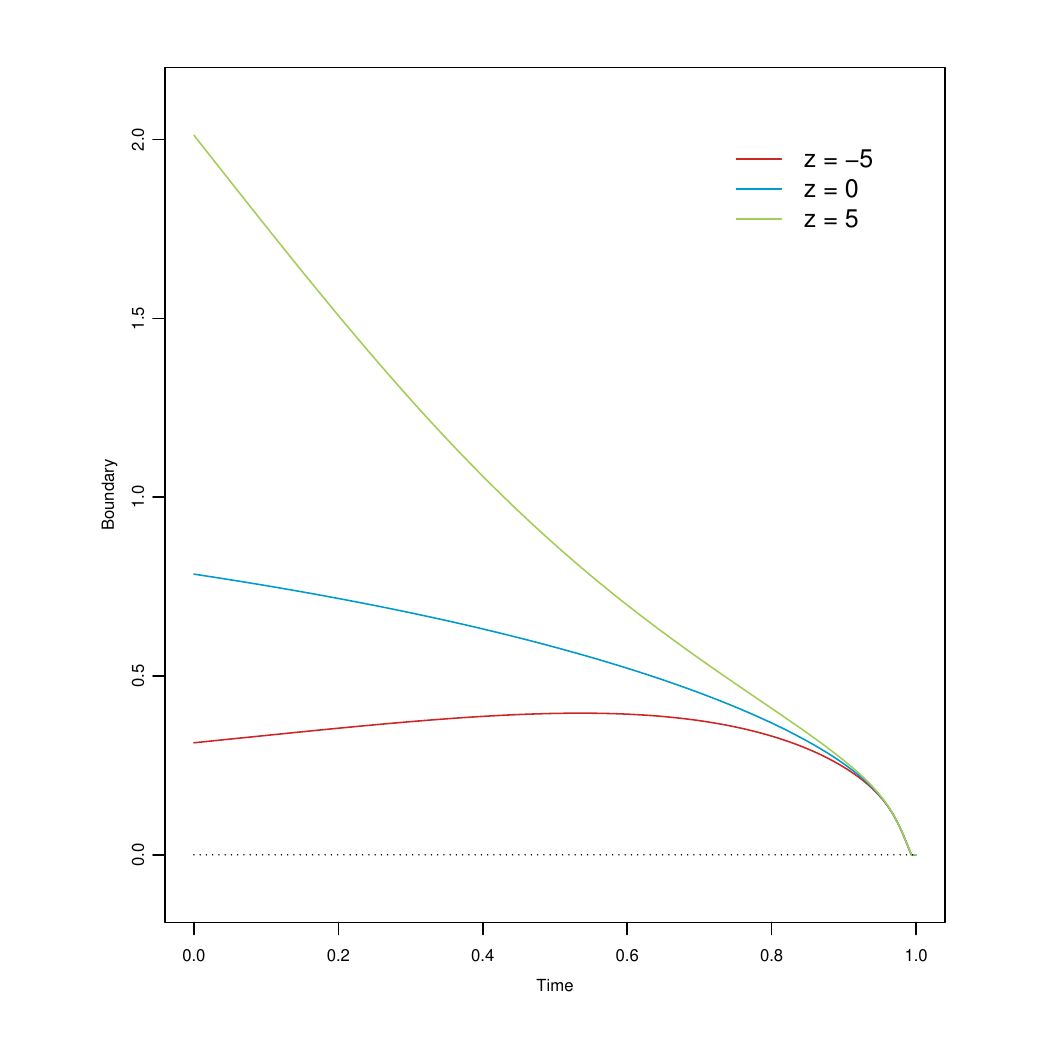}
		\subcaption{$\alpha = \gamma = 1, N = 100$}
	\end{subfigure}
	\begin{subfigure}[b]{0.32\textwidth}
		\includegraphics[width = \textwidth]{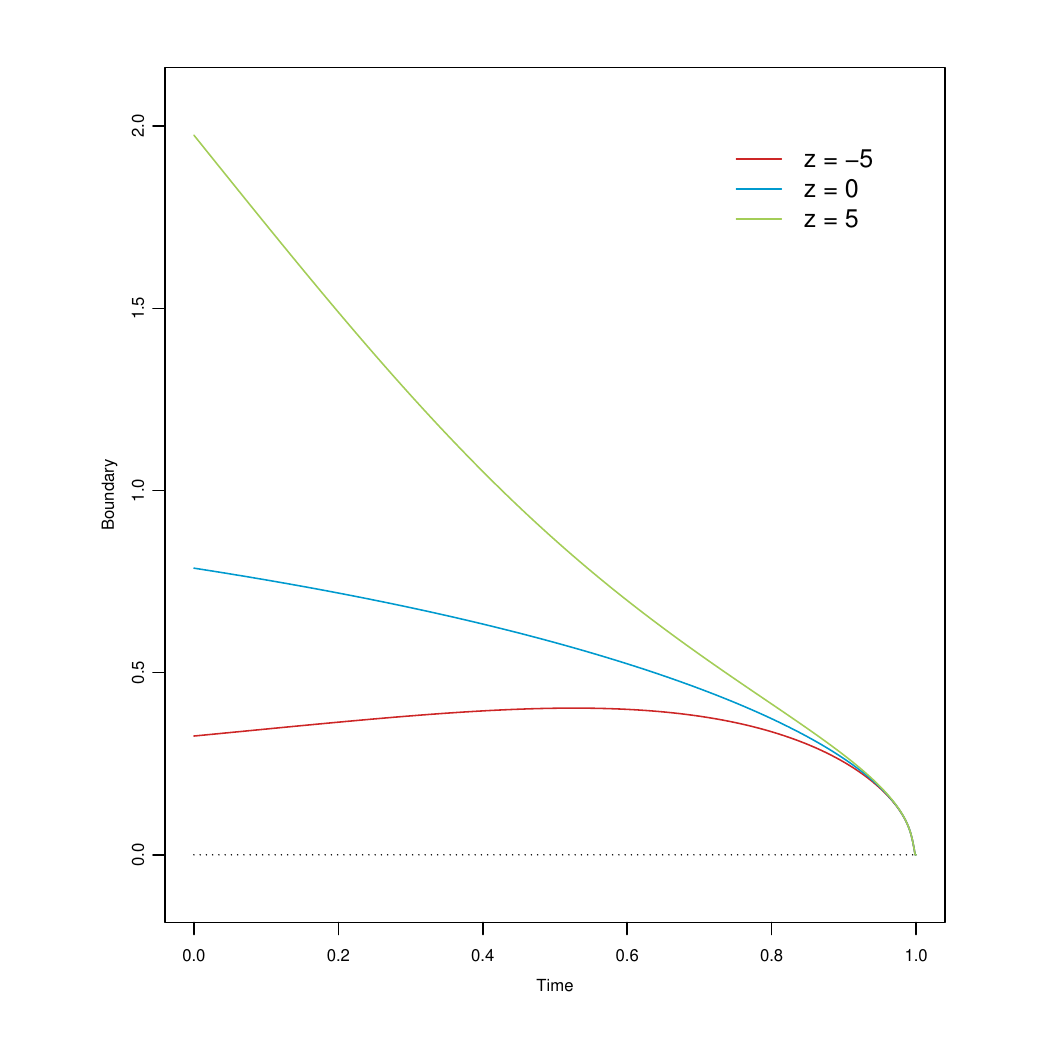}
		\subcaption{$\alpha = \gamma = 1, N = 500$}
	\end{subfigure}
	\vspace{6PT}

	\caption{\small Optimal stopping boundary estimation for different values of $z$ and $N$. We display $t\mapsto\beta(t) - z$ to allow a clearer comparison across the different values of $z$. As $N$ increases the boundary estimation is seen to converge.}
	\label{fig:pinning_change}
\end{figure}

The code implementing the boundary computation is available at \url{https://github.com/aguazz/OSP_OUB}.

%-------------------------------------------------%
\section{Conclusions}\label{sec:conclusions}
%-------------------------------------------------%

In this paper we solved the finite-horizon OSP for an OUB process with the identity as the gain function. To the best of our knowledge, so far the only Markov bridge addressed by the optimal stopping literature has been the BB and some slight variations of it (see, e.g., \cite{Shepp_1969_explicit, Follmer_1972_optimal, Ekstrom_2009_optimal, Ernst_2015_revisiting, Leung_2018_optimal, DeAngelis_2020_optimal, Glover_2020_optimally, Ekstrom_2020_optimal, DAuria_2020_discounted}). Markov bridges are potentially useful in mathematical finance as they allow including additional information at some terminal time.

Arguing as \cite{Shepp_1969_explicit} for the BB, we worked out the OUB case by coming up with an equivalent OSP having a Brownian motion as the underlying process after time-space transforming the OUB. Contrary to \cite{Shepp_1969_explicit}, the complexity of our problem did not allow us to guess a candidate solution, and we directly characterized the value function and the OSB by means of the pricing formula and the free-boundary equation. However, the equivalence between both OSPs was used only to facilitate technicalities along the proofs, and it is not necessary to compute the solution, since both the pricing formula and the free-boundary equation are also provided in the original formulation. We discussed how to use a Picard iteration algorithm to numerically approximate the OSB and displayed some examples to illustrate how different sets of values for the OUB's parameters rule the shape of the OSB.

%-------------------------------------------------%
\section*{Declaration of competing interest}
%-------------------------------------------------%
The authors declare that they have no known competing financial interests or personal relationships that could have appeared to influence the work reported in this paper.

%-------------------------------------------------%
\section*{Acknowledgments}
%-------------------------------------------------%

The authors acknowledge partial support from grants PID2020-116694GB-I00 (first and second authors), and PGC2018-097284-B-100 (third author), funded by MCIN/AEI/10.13039/\-501100011033 and by ``ERDF A way of making Europe''. The second author is a member of the Gruppo Nazionale Calcolo Scientifico-Istituto Nazionale di Alta Matematica (GNCS-INdAM) and acknowledges also the partial support from the Italian SID project BIRD239937/23.

\bibliographystyle{apalike-custom}
\bibliography{OUB-arXiv}

\end{document}